\documentclass[12pt]{article}
\usepackage{amsmath,amssymb,amsthm,euscript}

\makeatletter
\@addtoreset{equation}{section}

\newcommand{\gs}[1]{\textbf{#1}}
\newcommand{\ds}{\displaystyle}
\newcommand{\fr}[2]{\frac{#1}{#2}}
\newcommand{\dfr}[2]{\dfrac{#1}{#2}}
\newcommand{\cd}{\cdot}
\newcommand{\cds}{\cdots}
\newcommand{\dsum}{\displaystyle \sum}

\newcommand{\simto}{\stackrel{\sim}{\to}}

\renewcommand{\l}{\left}
\renewcommand{\r}{\right}

\newcommand{\vsv}{\vspace{5mm}}
\newcommand{\vsb}{\vspace{2mm}}

\newcommand{\q}{\quad}
\newcommand{\qq}{\qquad}

\newcommand{\maru}[1]{{\ooalign{\hfil#1\/\hfil\crcr
\raise.167ex\hbox{\mathhexbox20D}}}}

\newcommand{\ruby}[2]{%
 \leavevmode
 \setbox0=\hbox{#1}%
 \setbox1=\hbox{\tiny #2}%
 \ifdim\wd0>\wd1 \dimen0=\wd0 \else \dimen0=\wd1 \fi
 \hbox{%
   \kanjiskip=0pt plus 2fil
   \xkanjiskip=0pt plus 2fil
   \vbox{%
     \hbox to \dimen0{%
       \tiny \hfil#2\hfil}%
     \nointerlineskip
     \hbox to \dimen0{\mathstrut\hfil#1\hfil}}}}

\newcommand{\lfm}{\langle} 
\newcommand{\rfm}{\rangle} 
\newcommand{\la}{\langle}
\newcommand{\ra}{\rangle}

\newcommand{\abs}[1]{\lvert{#1}\rvert}

\DeclareMathOperator*{\tensor}{\otimes}
\DeclareMathOperator*{\fusion}{\boxtimes}

\newcommand{\Z}{\mathbb{Z}}
\newcommand{\C}{\mathbb{C}}
\newcommand{\R}{\mathbb{R}}
\newcommand{\N}{\mathbb{N}}
\newcommand{\Q}{\mathbb{Q}}

\newcommand{\aut}{\mathrm{Aut}}

\newcommand{\ind}{\mathrm{Ind}}

\newcommand{\w}{\omega}
\newcommand{\vacuum}{\mathrm{1\hspace{-3.2pt}l}}
\newcommand{\vac}{\vacuum}

\newcommand{\g}{\mathfrak{g}}

\newcommand{\hf}{\frac{1}{2}}

\newcommand{\pf}{\gs{Proof:}\q}
\newcommand{\M}{\mathbb{M}}

\newcommand{\mL}{\mathcal{L}}

\newcommand{\eh}{\mathbf{h}}
\newcommand{\ee}{\mathbf{e}}
\newcommand{\ef}{\mathbf{f}}
\newcommand{\nat}{\natural}

\newcommand{\Vef}{\mathrm{VA}(e,f)}

\theoremstyle{plain}
\newtheorem{thm}{Theorem}[section]
\newtheorem{prop}[thm]{Proposition}
\newtheorem{lem}[thm]{Lemma}

\theoremstyle{definition}
\newtheorem{df}[thm]{Definition}

\theoremstyle{remark}
\newtheorem{rem}[thm]{Remark}

\title{Vertex operator algebra with two Miyamoto involutions
generating $S_3$}
\date{}
\author{
  Shinya Sakuma
  \\
  \sf{\small sakuma@math.tsukuba.ac.jp}
  \vsb\\ 
  and 
  \vsb\\
  Hiroshi Yamauchi
  \\
  \sf{\small hirocci@math.tsukuba.ac.jp}
  \vsb\\
  \it{\small Graduate School of Mathematics,}
  \\
  \it{\small University of Tsukuba, Ibaraki 305-8571, Japan} 
}

\begin{document}

\baselineskip 6mm

\maketitle

\begin{abstract}
  \noindent
  In this article we study a VOA with two Miyamoto involutions 
  generating $S_3$. 
  In \cite{M3}, Miyamoto showed that a VOA generated by two conformal
  vectors whose Miyamoto involutions generate an automorphism group 
  isomorphic to $S_3$ is isomorphic to one of the four candidates he
  listed. 
  We construct one of them and prove that our VOA is actually the same 
  as $\mathrm{VA}(e,f)$ studied by Miyamoto.
  We also show that there is an embedding into the moonshine VOA.
  Using our VOA, we can define the 3A-triality in the Monster.
\end{abstract}

\section{Introduction} 

Vertex operator algebras (VOAs) associated with the unitary series of
the Virasoro algebras are very useful for studying VOA in which they
are contained.
This method was initiated by Dong et al. \cite{DMZ} in the study of
the moonshine VOA as a module of a tensor product of the first unitary 
Virasoro VOA $L(\hf,0)$.
Along this line, Miyamoto showed that the Virasoro VOA $L(\hf,0)$
defines an involution of a VOA in \cite{M1}, which is often called the
Miyamoto involution. 
In the moonshine VOA, this involution gives a 2A-involution in the
Monster sporadic simple group $\M$.
There are many interesting properties related to the 2A-involutions. 
For example, Mckay noted that there are some mysterious relations
between the $E_8$ Dynkin diagram and the 2A-involutions of the
Monster. There are also some similar relations between the
$Y_{555}$-diagram and the Bimonster.
For reference, see \cite{ATLAS} \cite{C} \cite{Mc}.
Motivated by the topics on 2A-involutions above, Miyamoto studied 
VOAs generated by two conformal vectors with central
charge $1/2$ whose Miyamoto involutions generate $S_3$ in \cite{M3}
and he determined that the possible inner products of such a pair of
conformal vectors are $1/2^8$ or $13/2^{10}$.
Furthermore, he determined the possible candidates of VOAs
generated by such two conformal vectors.
When the inner product is equal to $13/2^{10}$, he showed 
that a VOA generated by such two conformal vectors is isomorphic to
one of the following (Theorem 5.6 of \cite{M3}):
\vsb\\
(1)\ 
  $
    \l( L(\fr{4}{5},0)\oplus L(\fr{4}{5},3)\r) \tensor L(\fr{6}{7},0) 
    \oplus L(\fr{4}{5},\fr{2}{3})^+\tensor L(\fr{6}{7},\fr{4}{3})
    \oplus L(\fr{4}{5},\fr{2}{3})^-\tensor L(\fr{6}{7},\fr{4}{3}),
  $
  \vsb\\
(2)\ 
  $
    L(\fr{4}{5},0)\tensor \l( L(\fr{6}{7},0)\oplus L(\fr{6}{7},5)\r)
    \oplus L(\fr{4}{5},\fr{2}{3})\tensor L(\fr{6}{7},\fr{4}{3})^+
    \oplus L(\fr{4}{5},\fr{2}{3})\tensor L(\fr{6}{7},\fr{4}{3})^-,
  $
  \vsb\\
(3)\
  $
    L(\fr{4}{5},0)\tensor L(\fr{6}{7},0)
    \oplus L(\fr{4}{5},3)\tensor L(\fr{6}{7},5)
    \oplus \l( L(\fr{4}{5},\fr{2}{3})\tensor L(\fr{6}{7},\fr{4}{3})
      \r)^+ 
    \oplus \l( L(\fr{4}{5},\fr{2}{3})\tensor L(\fr{6}{7},\fr{4}{3})
      \r)^-,
  $
  \vsb\\
(4)\ 
  $
    \l( L(\fr{4}{5},0)\oplus L(\fr{4}{5},3)\r)
    \tensor \l( L(\fr{6}{7},0)\oplus L(\fr{6}{7},5)\r)
    \oplus L(\fr{4}{5},\fr{2}{3})^+ \tensor L(\fr{6}{7},\fr{4}{3})^\pm
    \oplus L(\fr{4}{5},\fr{2}{3})^-\tensor L(\fr{6}{7},\fr{4}{3})^\mp .
  $
\vsb\\
Unfortunately, these VOAs are just candidates and it is still
unknown if they actually exist.
In this paper, we construct a VOA $U$ which has the same shape as that
of the candidate (4).
We show that in (4) there is a unique simple VOA structure and 
$(++)$-type structure and $(+-)$-type structure are isomorphic and 
conjugate of each other. 
We classify all irreducible modules and the fuion algebra for $U$
and prove that it is a rational VOA.
We also prove that it is generated by two conformal vectors with
central charge $1/2$ whose inner product is $13/2^{10}$ and also we
show that their Miyamoto involutions generate $S_3$. 
Namely, $U$ is the same as the VOA studied in \cite{M3} and gives 
a positive solution for Theorem 5.6 (4) of \cite{M3}.
We further prove that the candidates (1)-(3) do not exist so that
only the candidate (4) occurs (Theorem \ref{completion}).
Therefore, we can verify that $U$ is contained in the moonshine VOA
and using fusion rules for $U$-modules we can define the 3A-triality
of the Monster (Theorem \ref{3A}).
Throughout the paper, we will work over the field $\C$ of complex
numbers unless otherwise stated.

\section{Preliminaries} 

\subsection{The unitary series of the Virasoro VOAs} 

For any complex numbers $c$ and $h$, denote by $L(c,h)$ the
irreducible highest weight representation of the Virasoro algebra with
central charge $c$ and highest weight $h$.
It is shown in \cite{FZ} that $L(c,0)$ has a natural structure of a
simple VOA.
Let 
\begin{align}
  c_m &:= 1-\dfr{6}{(m+2)(m+3)}\qq (m=1,2,\dots ),
  \vsb\\
  h_{r,s}^{(m)} &:= \dfr{\{ r(m+3)-s(m+2)\}^2-1}{4(m+2)(m+3)} 
  \label{highest weight}
\end{align}
for $r,s\in \N$, $1\leq r\leq m+1$ and $1\leq s\leq m+2$.
It is shown in \cite{W} that $L(c_m,0)$ is rational and 
$L(c_m,h_{r,s}^{(m)})$, $1\leq s\leq r\leq m+1$, provide all
irreducible $L(c_m,0)$-modules (see also \cite{DMZ}).
This is so-called the unitary series of the Virasoro VOAs.
The fusion rules among $L(c_m,0)$-modules \cite{W} are given by 
\begin{equation}\label{wang}
  L(c_m,h_{r_1,s_1})\times L(c_m,h_{r_2,s_2}) 
  = \dsum_{i\in I,j\in J} L(c_m,h_{\abs{r_1-r_2}+2i-1,\abs{s_1-s_2}
  +2j-1}), 
\end{equation}
where 
$$
\begin{array}{l}
  I=\{ 1,2,\dots,\min \{ r_1,r_2,m+2-r_1,m+2-r_2\}\} ,
  \vsb\\
  J=\{ 1,2,\dots,\min \{ s_1,s_2,m+3-s_1,m+3-s_2\}\} .
\end{array}
$$

\subsection{GKO-construction} 

Let $\g$ be the Lie algebra $\hat{\textit{sl}_2}(\C )$ with 
generators $h$, $e$, $f$ and relations $[h,e]=2e$, $[h,f]=-2f$ and
$[e,f]=h$. 
We use the standard invariant bilinear form on $\g$ defined by
$\la h,h\ra =2$ and $\la e,f\ra =1$.
Let $\hat{\g}$ be the corresponding affine algebra of type $A_1^{(1)}$ 
and $\Lambda_0$, $\Lambda_1$ the fundamental weights for $\hat{\g}$.
For any non-negative integers $m$ and $j$, denote by
$\mathcal{L}(m,j)$ the irreducible highest weight $\hat{\g}$-module
with highest weight $(m-j)\Lambda_0 +j\Lambda_1$.
Then $\mL (m,0)$ has a natural structure of a simple VOA \cite{FZ}.
The Virasoro vector $\Omega^m$ of $\mL (m,0)$ is given by
\begin{equation}\label{Omega}
  \Omega^m := \dfr{1}{2(m+2)}\l( \fr{1}{2}h_{(-1)} h + e_{(-1)} f
  +f_{(-1)} e \r)  
\end{equation}
with central charge $3m/(m+2)$.

Let $m\in \N$. Then $\mL (m,0)$ is a rational VOA and  
$\{ \mL (m,j) \mid j=0,1,\dots,m\}$ is the set of all irreducible 
$\mL (m,0)$-modules.  
The fusion algebra (cf. \cite{FZ}) is given by
\begin{equation}\label{FZ}
  \mL (m,j) \times \mL (m,k) 
  = \sum_{i=\max \{ 0,j+k-m\}}^{\min \{j,k\}} \mL (m,j+k-2i).
\end{equation}
In particular, $\mL (m,m)\times \mL (m,j) = \mL (m,m-j)$ and thus
$\mL (m,m)$ is a simple current module.
A reasonable explanation why $\mL (m,m)$ is a simple current is given
in \cite{Li2}.

The weight 1 subspace of $\mL (m,0)$ forms a Lie algebra isomorphic to 
$\g$ under the $0$-th product in $\mL (m,0)$.
Let $\eh^1,\ee^1,\ef^1$ be the generator of $\g$ in $\mL (1,0)_1$ and 
$\eh^m,\ee^m,\ef^m$ those in $\mL (m,0)_1$.
Then  $\eh^{m+1}:= \eh^1\tensor \vac +\vac \tensor \eh^m$,
$\ee^{m+1}:=\ee^1\tensor \vac + \vac \tensor \ee^m$ and $\ef^{m+1}:=
\ef^1\tensor \vac + \vac \tensor \ef^m$ generate a sub VOA isomorphic
to $\mL (m+1,0)$ in $\mL (1,0)\tensor \mL (m,0)$ with the Virasoro
vector $\Omega^{m+1}$ made from $\eh^{m+1}$, $\ee^{m+1}$ and
$\ef^{m+1}$ by \eqref{Omega}. 
It is shown in \cite{DL} and \cite{KR} that $\w^m:= \Omega^1\tensor
\vac + \vac \tensor \Omega^m -\Omega^{m+1}$ also gives a Virasoro
vector with central charge $c_m=1-6/(m+2)(m+3)$.
Furthermore, $\Omega^{m+1}$ and $\w^m$ are mutually commutative and
$\w^m$ generates a simple Virasoro VOA $L(c_m,0)$. 
Hence, $\mL (1,0)\tensor \mL (m,0)$ contains a sub VOA isomorphic to 
$L(c_m,0)\tensor \mL (m+1,0)$.
Since both $L(c_m,0)$ and $\mL (m+1,0)$ are rational,
every $\mL (1,0)\tensor \mL (m,0)$-module can be decomposed into
irreducible $L(c_m,0)\tensor \mL (m+1,0)$-submodules.
The  following decomposition is obtained in \cite{GKO}:
\begin{equation}\label{character}
  \mL (1,\epsilon )\tensor \mL (m,n) 
  = \bigoplus_{0\leq s\leq m+1 \atop s\equiv n+\epsilon \mod 2}
    L(c_m,h_{n+1,s+1}^{(m)})\tensor \mL (m+1,s),
\end{equation}
where $\epsilon =0,1$ and $0\leq n\leq m$.
Note that $h_{r,s}^{(m)}=h_{m+2-r,m+3-s}^{(m)}$.
This is the famous GKO-construction of the unitary Virasoro VOAs.

\subsection{Lattice construction of $\mL (m,0)$} 

Let $A_1=\Z \alpha$ with $\lfm \alpha,\alpha\rfm =2$ be the root
lattice of type $A_1$ and $V_{A_1}$ the lattice VOA associated with
$A_1$. 
Let
$$
  A_1^*=\{ x\in \Q \tensor_\Z A_1 \mid \lfm x,\alpha \rfm \in \Z\}
$$
be the dual lattice of $A_1$.
Then $A_1^*=A_1\cup (\hf \alpha +A_1)$.
It is well-known that $V_{A_1}\simeq \mL (1,0)$ and $V_{\hf \alpha
+A_1}\simeq \mL (1,1)$ (cf. \cite{FLM} \cite{FZ}, etc.).
Let $A_1^m=\Z \alpha^1\oplus \Z \alpha^2 \oplus \cds \oplus \Z
\alpha^m$ be the orthogonal sum of $m$ copies of $A_1$. 
Then we have an isomorphism $V_{A_1^m} \simeq (V_{A_1})^{\tensor m}
\simeq \mL (1,0)^{\tensor m}$.
Let 
$H^m:= \alpha^1_{(-1)}\vac +\cds + \alpha^m_{(-1)}\vac$, 
$E^m:= e^{\alpha^1}+\cds + e^{\alpha^m}$ and
$F^m:= e^{-\alpha^1}+\cds + e^{-\alpha^m}$.
Then it is shown in \cite{DL} that $H^m$, $E^m$ and $F^m$ generate 
a sub VOA isomorphic to $\mL (m,0)$ in $V_{A_1^m}$.

\subsection{Vertex operator algebra $L(\fr{4}{5},0)\oplus 
L(\fr{4}{5},3)$} 

Here we review the simple VOA $L(\fr{4}{5},0)\oplus
L(\fr{4}{5},3)$.
It is a $\Z_2$-simple current extension of the unitary Virasoro VOA
$L(\fr{4}{5},0)$ and is deeply studied in \cite{KMY} and \cite{M2}.
By the fusion rule \eqref{wang}, there exists a canonical involution
$\sigma$ on $L(\fr{4}{5},0)\oplus L(\fr{4}{5},3)$ which acts as
identity on $L(\fr{4}{5},0)$ and acts as a scalar $-1$ on
$L(\fr{4}{5},3)$. 
We also note that $\sigma$ is the only non-trivial automorphism on 
$L(\fr{4}{5},0)\oplus L(\fr{4}{5},3)$.
For any $L(\fr{4}{5},0)\oplus L(\fr{4}{5},3)$-module $(M,Y_M(\cd,z))$, 
we can consider its {\it $\sigma$-conjugate module} $(M^\sigma,
Y_M(\cd,z))$ which is defined as follows.
As a vector space, we put $M^\sigma \simeq M$ and the action of $a\in
L(\fr{4}{5},0)\oplus L(\fr{4}{5},3)$ is given by 
$$
  Y_M^\sigma (a,z):= Y_M(\sigma a,z).
$$
We will denote the $\sigma$-conjugate of $M$ simply by $M^\sigma$.

\begin{thm}
  (\cite{KMY})\
  A VOA $L(\fr{4}{5},0)\oplus L(\fr{4}{5},3)$ is rational and every 
  irreducible module is isomorphic to one of the following:
  $$
  \begin{array}{ll}
    W(0):= L(\fr{4}{5},0)\oplus L(\fr{4}{5},3),
    & W(\fr{2}{3})^\pm := L(\fr{4}{5},\fr{2}{3})^\pm,
    \vsb\\
    W(\fr{2}{5}):= L(\fr{4}{5},\fr{2}{5})\oplus
      L(\fr{4}{5},\fr{7}{5}),    
    & W(\fr{1}{15})^\pm := L(\fr{4}{5},\fr{1}{15})^\pm,
  \end{array}
  $$
  where $W(h)^-$ is the $\sigma$-conjugate module of $W(h)^+$.
  The dual modules are as follows:
  $(W(h)^\pm)^*\simeq W(h)^\mp$ if $h=\fr{2}{3}$ or $\fr{1}{15}$ 
  and $W(h)^*\simeq W(h)$ for the others.
\end{thm}

\begin{rem}\label{signature}
  We may exchange the sign $\pm$ since there is no canonical way
  to determine the type $+$ and $-$ for the modules $W(h)^+$ and
  $W(h)^-$. 
  However, if we determine a sign of one module, then the
  following fusion rules automatically determine all the signs.
\end{rem}

The fusion algebra for $W(0)$ has a natural $\Z_3$-symmetry.
For convenience, we use the following $\Z_3$-graded names.
$$
\begin{array}{lll}
  A^0:= W(0), & A^1:= W(\fr{2}{3})^+, & A^2:= W(\fr{2}{3})^-,
  \vsb\\
  B^0:= W(\fr{2}{5}), & B^1:= W(\fr{1}{15})^+, & B^2:=
    W(\fr{1}{15})^-.
\end{array}
$$

\begin{thm}\label{fusion for W(0)}
  (\cite{M2})\
  The fusion rules for irreducible $W(0)$-modules are given as
  $$
  \begin{array}{l}
    A^i \times A^j = A^{i+j},
    \vsb\\
    A^i \times B^j = B^{i+j},
    \vsb\\
    B^i \times B^j = A^{i+j} + B^{i+j},
  \end{array}
  $$
  where $i,j\in \Z_3$. Therefore, the fusion algebra for $W(0)$ has a
  natural $\Z_3$-symmetry.
\end{thm}

\subsection{Vertex operator algebra $L(\fr{6}{7},0) \oplus 
L(\fr{6}{7},5)$} 

In this subsection we give some facts about the VOA
$L(\fr{6}{7},0)\oplus L(\fr{6}{7},5)$.
This is a $\Z_2$-simple current extension of the unitary Virasoro VOA
$L(\fr{6}{7},0)$ and is studied in \cite{LY}.
Also, all statements in this subsection are included in
\cite{LLY}.
So we give a slight explanation here.

\begin{thm}\label{tricritical}
  (\cite{LY} and \cite{LLY})\ 
  There exists a unique structure of a simple VOA on $L(\fr67,0)\oplus 
  L(\fr67,5)$.
\end{thm}

\pf
It follows from the fusion rules \eqref{wang} that if it has a
structure of a VOA then it must be unique.
So we should show the existence of a structure.
This will be given later.
\qed

As in the case of $L(\fr{4}{5},0)\oplus L(\fr{4}{5},3)$, 
a linear map $\sigma$ which acts as a scalar 1 on $L(\fr{6}{7},0)$ and 
acts as $-1$ on $L(\fr{6}{7},5)$ defines an automorphism of a VOA
$L(\fr{6}{7},0)\oplus L(\fr{6}{7},5)$. 
We also note that $\sigma$ is the only non-trivial automorphism on
$L(\fr{6}{7},0)\oplus L(\fr{6}{7},5)$.

\begin{thm}(\cite{LY} and \cite{LLY})\ 
  A VOA $L(\fr{6}{7},0)\oplus L(\fr{6}{7},5)$ is rational and 
  all its irreducible modules are the following:
  $$
  \begin{array}{lll}
    N(0):=L(\fr{6}{7},0)\oplus L(\fr{6}{7},5),
    & N(\fr{1}{7}):= L(\fr{6}{7},\fr{1}{7})\oplus
      L(\fr{6}{7},\fr{22}{7}) ,
    & N(\fr{5}{7}):= L(\fr{6}{7},\fr{5}{7})\oplus
      L(\fr{6}{7},\fr{12}{7}), 
    \vsb\\
    N(\fr{4}{3})^\pm := L(\fr{6}{7},\fr{4}{3})^\pm,
    & N(\fr{1}{21})^\pm := L(\fr{6}{7},\fr{1}{21})^\pm,
    & N(\fr{10}{21})^\pm := L(\fr{6}{7},\fr{10}{21})^\pm ,
  \end{array}
  $$
  where $N(h)^-$ is the $\sigma$-conjugate module of $N(h)^+$.
  Also, the dual modules are as follows:
  $(N(h)^\pm )^*\simeq N(h)^\mp$ if $h= \fr{4}{3}$, $\fr{1}{21}$ or
  $\fr{10}{21}$ and $N(h)^*\simeq N(h)$ for the others. 
\end{thm}

The fusion algebra for $N(0)$ is also determined in \cite{LY}
and \cite{LLY}.
To state the fusion rules, we assign $\Z_3$-graded names to
irreducible modules (cf. \cite{LY}). 
Define
$$
\begin{array}{lll}
  A^0 := N(0), & A^1 := N(\fr{4}{3})^+, & A^2 :=N(\fr{4}{3})^-,
  \vsb\\
  B^0 := N(\fr{1}{7}), & B^1 := N(\fr{10}{21})^+, 
  & B^2 := N(\fr{10}{21})^-, 
  \vsb\\
  C^0 := N(\fr{5}{7}), & C^1 := N(\fr{1}{21})^+, 
  & C^2 := N(\fr{1}{21})^-. 
\end{array}
$$

\begin{thm} \label{fusion for N(0)}
  (\cite{LY} and \cite{LLY})\ 
  The fusion rules for irreducible $N(0)$-modules are given as 
  $$
  \begin{array}{lll}
    A^i \times A^j &=& A^{i+j},
    \vsb\\
    A^i \times B^j &=& B^{i+j},
    \vsb\\
    A^i \times C^j &=& C^{i+j},
    \vsb\\
    B^i \times B^j &=& A^{i+j} + C^{i+j},
    \vsb\\
    B^i \times C^j &=& B^{i+j} + C^{i+j},
    \vsb\\
    C^i \times C^j &=& A^{i+j} + B^{i+j} + C^{i+j},
  \end{array}
  $$
  where $i,j\in \Z_3$. Therefore, the fusion algebra for $N(0)$ has a
  natural $\Z_3$-symmetry.
\end{thm}

\section{Simple current extensions} 

In this section we consider how vertex operator algebras are
extended by their simple current modules (see also \cite{L}).
Let $D$ be an abelian group and $V^0$ a simple and rational VOA.
Assume that a set of irreducible $V^0$-modules $\{ V^\alpha \mid
\alpha \in D\}$ indexed by $D$ is given.
One can easily verify the following lemma.

\begin{lem}
  Assume that $\oplus_{\alpha\in D}V^\alpha$ carries a structure of a 
  $D$-graded VOA such that $V^\alpha \cd V^\beta =\{ \sum a_{(n)} b 
  \mid a\in V^\alpha, b\in V^\beta, n\in \Z \} \ne 0$.
  It is simple if and only if $V^\alpha$ and $V^\beta$ are 
  inequivalent irreducible $V^0$-modules for distinct $\alpha$ and 
  $\beta\in D$.
\end{lem}

\pf
Assume that $V_D$ is simple.
Then the automorphism group of $V_D$ contains a group isomorphic to
the dual group $D^*$ of an abelian group $D$ because $V_D$ is
$D$-graded.
It is clear that the $D^*$-invariants of $V_D$ is exactly $V^0$.
Therefore, by the quantum Galois theory \cite{DM1} \cite{HMT}, each
$V^\alpha$ is an irreducible $V^0$-modules.

Conversely, if $\{ V^\alpha \mid \alpha \in D\}$ is a set of
inequivalent irreducible $V^0$-modules such that
$V_D=\oplus_{\alpha\in D}V^\alpha$ forms a $D$-graded vertex operator
algebra, then $V_D$ must be simple because of the density theorem.
\qed

The lemma above leads us the following definition.
\begin{df}
  A {\it $D$-graded extension} $V_D$ of $V^0$ is a simple VOA with the
  shape $V_D=\oplus_{\alpha\in D}V^\alpha$ whose vacuum element and
  Virasoro element are given by those of $V^0$ and vertex operations
  in $V_D$ satisfies $Y(u^\alpha,z)v^\beta \in V^{\alpha +\beta}((z))$
  for any $u^\alpha\in V^\alpha$ and $v^\beta \in V^\beta$.
\end{df}

If a $D$-graded extension $V_D$ of $V^0$ is given,
then the rationality of $V^0$ implies that $D$ is a finite abelian.
It is natural for us to ask how many structures can sit in $V_D$.

\begin{lem}\label{uniqueness0}
  (\cite[Proposition 5.3]{DM2})
  Suppose that the space of $V^0$-intertwining operators of type
  $V^\alpha\times V^\beta \to V^{\alpha+\beta}$ is one dimensional.
  Then the VOA structure of a $D$-graded extension $V_D$ of $V^0$ over 
  $\C$ is unique.
\end{lem}

\pf
Assume that two $D$-graded extensions
$V_D=(\oplus_{\alpha \in D} V^\alpha ,Y^1(\cd,z))$ and
$\tilde{V}_D=(\oplus_{\alpha\in D} \tilde{V}^\alpha
,Y^2(\cd,z))$ of $V^0$ such that $V^\alpha \simeq \tilde{V}^\alpha$ 
as $V^0$-modules are given.
Then there exist $V^0$-isomorphisms $\phi_\alpha : V^\alpha \to
\tilde{V}^\alpha$ such that $Y^2(\phi_0 a,z)\phi_\alpha = \phi_\alpha
Y^1(a,z)$ for all $a\in V^0$.
Since both $Y^1(\cd,z)|_{V^\alpha \tensor V^\beta}$ and 
$\phi_{\alpha+\beta}^{-1}
Y^2(\phi_\alpha\, \cd,z)\phi_\beta |_{V^\alpha\tensor V^\beta}$ are
$V^0$-intertwining operators of type $V^\alpha \times V^\beta \to
V^{\alpha +\beta}$, there exist scalars $c(\alpha,\beta)\in \C$ such
that  
\begin{equation}
  Y^2(\phi_\alpha a,z) \phi_\beta b 
  = c(\alpha,\beta)\phi_{\alpha +\beta} Y^1(a,z)b
\end{equation}
for all $a\in V^\alpha$ and $b\in V^\beta$.
Recall that the associativity and the commutativity of vertex
operators.
Let $x,y$ be any element in a VOA and $v$ be any element in a module.
Then there exist $N_1, N_2\in \N$ such that 
\begin{align}
  & (z_1-z_2)^{N_1} Y(x,z_1)Y(y,z_2) v
    = (z_1-z_2)^{N_1} Y(x,z_2) Y(y,z_1)v, 
    \label{commutativity}
  \vsb\\
  & (z_0+z_2)^{N_2} Y(x,z_0+z_2)Y(y,z_2)v
    = (z_2+z_0)^{N_2} Y(Y(x,z_0)y,z_2)v.
    \label{associativity}
\end{align}
The first equality is called the commutativity and the second is
called the associativity of vertex operators.
An integer $N_1$ depends on $x$ and $y$, whereas $N_2$ does not only 
on $x$ and $y$ but also $v$.
Using the commutativity \eqref{commutativity} and the associativity 
\eqref{associativity}, we can show that $c(\cd,\cd): D\times D \to \C$ 
satisfies the following:
\begin{equation}
  c(\alpha,\beta+\gamma)c(\beta,\gamma)
  = c(\beta,\alpha+\gamma)c(\alpha,\gamma)
  = c(\alpha,\beta)c(\alpha+\beta,\gamma).
\end{equation}
Moreover, by the definition and by the skew-symmetry, we also have
$c(0,\alpha)=1$ and $c(\alpha,\beta)=c(\beta,\alpha)$ for all
$\alpha,\beta\in D$.
Namely, $c(\cd,\cd)$ is a 2-cocycle on $D$ and hence it defines a
central extension of $D$ by $\C$.
Since $D$ is a direct sum of finite cyclic groups and $\C$ is
algebraically closed, such a central extension splits and hence there
exists a coboundary $t: D\to \C$ such that $c(\alpha,\beta) =
t(\alpha+\beta) t(\alpha)^{-1} t(\beta)^{-1}$ for all $\alpha,\beta\in
D$. 
Then by putting $\tilde{\phi}|_{V^\alpha}:= t(\alpha)\phi_\alpha$, 
we obtain a VOA-isomorphism $\tilde{\phi} : V_D \simto \tilde{V}_D$. 
This completes the proof.
\qed

By the lemma above, we adopt the following definitions.
\begin{df}
  An irreducible $V^0$-module $X$ is called a {\it simple current
  $V^0$-module} if it satisfies that for every irreducible
  $V^0$-module $W$, the fusion product (or the tensor product) 
  $X\times W$ is also irreducible.
\end{df}

\begin{df}
  A $D$-graded extension $V_D=\oplus_{\alpha\in D} V^\alpha$ of $V^0$
  is called a $D$-graded {\it simple current extension} if all
  $V^\alpha$, $\alpha\in D$, are simple current $V^0$-modules.
\end{df}

Clearly, if $V_D$ is a $D$-graded simple current extension, then it
satisfies the assumption in Lemma \ref{uniqueness0}.
Let $E$ be any subgroup of $D$ and $D=\bigcup_{i=1}^{\abs{D}/\abs{E}}
(t^i+E)$ a coset decomposition of $D$ with respect to $E$.
Set $V^{t_i+E}:= \oplus_{\beta\in E} V^{t^i+\beta}$.
The definition of $V^{t_i+E}$ does not depend on the choice of
representatives $\{ t_i\}$.
It is clear from the definition that $V_E:=\oplus_{\alpha\in
E}V^\alpha$ is an $E$-graded extension of $V^0$ and
$V_{D/E}:=\oplus_{i=1}^{\abs{D}/\abs{E}} V^{t_i+E}$ is a $D/E$-graded
extension of $V_E$.
Furthermore, if $V_D$ is a $D$-graded simple current extension, then 
$V_E$ (resp. $V_{D/E}$) is also an $E$-graded and (resp. $D/E$-graded)
simple current extension of $V^0$ (resp. $V_E$); the proof will be
given in Lemma \ref{restriction}.
See Remark \ref{SCE}.

Let $M$ be an irreducible $V_D$-module.
Since we have assumed that $V^0$ is rational, there is an irreducible
$V^0$-submodule $W$ of $M$.

\begin{lem}\label{non-zero}
  Let $V_D$ be a $D$-graded extension of $V^0$ and let $M$ be an 
  admissible $V_D$-module. 
  For an irreducible $V^0$-submodule $W$ of $M$, 
  $V^\alpha \cd W:=\{ \sum a_{(n)} w \mid a\in V^\alpha, w\in W, n\in 
  \Z\}$ are also non-trivial irreducible $V^0$-submodules for all
  $\alpha\in D$. 
\end{lem}

\pf
First, note that $V^\alpha\cd W$ is a $V^0$-submodule by
\eqref{associativity}.
Moreover, the associativity \eqref{associativity} tells us that
$V^\alpha \cd (V^\beta \cd W)\subset (V^\alpha \cd V^\beta)\cd W=
V^{\alpha +\beta}\cd W$.
We show that $V^\alpha\cd W$ is not zero and then we prove that it is
irreducible.
If $V^\alpha \cd W=0$, then by the iterate formula
$$
  (a_{(m)}b)_{(n)} = \sum_{i=0}^\infty (-1)^i \binom{m}{i} 
  \{ a_{(m-i)} b_{(n+i)}-(-1)^m b_{(m+n-i)}a_{(i)}\}
$$
we obtain $V^{n\alpha}\cd W=0$ for $n=1,2,\dots$. 
But $D$ is a finite abelian, we arrive at $V^0\cd W=0$, a
contradiction. 
Therefore, $V^\alpha\cd W\ne 0$ for all $\alpha \in D$.
Next, assume that there exists a proper non-trivial $V^0$-submodule
$X$ in  $V^\alpha \cd W$. 
Then we have $V^{-\alpha}\cd X \subset V^{-\alpha}\cd (V^\alpha \cd W) 
\subset (V^{-\alpha} \cd V^{\alpha})\cd W = V^0\cd W=W$ and hence we
get $V^{-\alpha}\cd X=W$ because $W$ is irreducible.
Then we obtain $V^\alpha \cd W = V^\alpha \cd (V^{-\alpha}\cd X)\subset
(V^\alpha \cd V^{-\alpha})\cd X= V^0\cd X=X$, a contradiction.
Therefore, $V^\alpha\cd W$ is a non-trivial and irreducible
$V^0$-submodule of $M$.
\qed

Let $M$ and $W$ be as in the lemma above and assume that $M$ is
irreducible under $V_D$.
Then $M=V_D\cd W=\sum_{\alpha \in D} V^\alpha \cd W$.
Set $D_W:= \{\alpha \in D \mid V^\alpha\cd W\simeq W\}$.
Since both $V^\alpha\cd (V^\beta\cd W)$ and $V^{\alpha +\beta}\cd W$
are irreducible $V^0$-modules by the previous lemma, it follows from
the associativity that $D_W$ is a subgroup of $D$.
Let $D=\bigcup_{i=1}^{\abs{D/D_W}} (\alpha^i+D_W)$ be a coset
decomposition with $\alpha^1=0$.
We note that $V^\alpha \cd W \simeq V^\beta \cd W$ if and only if 
$\alpha \in \beta +D_W$.
Set $M^{\alpha^i+D_W}:= \sum_{\beta\in D_W} (V^{\alpha^i+\beta}\cd W)$.
Then $M^{\alpha^i+D_W}$ is a direct sum of some copies of
$V^{\alpha^i}\cd W$'s as a $V^0$-module and $M$ decomposes into a
direct sum of $\abs{D/D_W}$-isotypical components 
$$
  M=\bigoplus_{i=1}^{\abs{D/D_W}} M^{\alpha^i+D_W}
$$
as a $V^0$-module.
We note that each $M^{\alpha^i+D_W}$ is a $V_{D_W}$-module and 
$M$ is a $D/D_W$-graded $V_{D_W}$-module, that is,
$V^{\alpha^i+D_W}\cd M^{\alpha^j+D_W}= M^{\alpha^i+\alpha^j+D_W}$.
Therefore, by the irreducibility of $M$, all $M^{\alpha^i+D_W}$ are
irreducible $V_{D_W}$-submodules.

\begin{df}
  A $V_D$-module $M$ is said to be {\it $D$-stable} if $D_W=0$ for
  some irreducible $V^0$-submodule $W$ of $M$.
\end{df}

It is obvious that the definition of the $D$-stability is independent
of the choice of an irreducible $V^0$-module $W$.

\begin{prop}\label{unique1}
  Let $V_D$ be a $D$-graded simple current extension of $V^0$. 
  Then the structure of every $D$-stable irreducible $V_D$-module is
  unique over $\C$.
  In other words, the $V^0$-module structure completely determines 
  the $V_D$-module structure of all $D$-stable irreducible
  $V_D$-modules. 
\end{prop}

\pf
Let $M$ be a $D$-stable irreducible $V_D$-module and let $W$ be an
irreducible $V^0$-submodule of $M$.
By definition, we have $M=\oplus_{\alpha \in D} (V^\alpha \cd W)$ and
all $V^\alpha\cd W$, $\alpha\in D$, are non-trivial irreducible
$V^0$-submodules. 
Set $W^\alpha:= V^\alpha \cd W$ for $\alpha\in D$.
We show that there exists a unique $V_D$-module structure on
$\oplus_{\alpha\in D}W^\alpha$.
Suppose that there are two $V_D$-modules $M=(\oplus_{\alpha\in D}
W^\alpha, Y^1(\cd,z))$ and $\tilde{M}=(\oplus_{\alpha\in D}
\tilde{W}^\alpha, Y^2(\cd,z))$ such that $W^\alpha\simeq
\tilde{W}^\alpha$ as $V^0$-modules for all $\alpha\in D$.
By assumption, there exist $V^0$-isomorphism $\psi_\alpha :
W^\alpha\to \tilde{W}^\alpha$ such that $Y^2(a,z)\psi_\alpha =
\psi_\alpha Y^1(a,z)$ for all $a\in V^0$.
Then both $Y^1(\cd,z)|_{V^\alpha\tensor W^\beta}$ and $\psi_{\alpha
+\beta}^{-1} Y^2(\cd ,z)\psi_\beta |_{V^\alpha \tensor
W^\beta}$ are $V^0$-intertwining operators of type
$V^\alpha \times W^\beta \to W^{\alpha +\beta}$ and hence there exist
non-zero scalars $c(\alpha,\beta)\in \C$ such that
$Y^2(a,z)\psi_\beta = c(\alpha,\beta) \psi_{\alpha+\beta} Y^1(a,z)$
for all $a\in V^\alpha$.
Then, by the associativity \eqref{associativity} we obtain
\begin{equation}\label{cocycle}
  c(\alpha +\beta,\gamma)
  = c(\alpha,\beta +\gamma)c(\beta,\gamma)
\end{equation}
for $\alpha, \beta, \gamma\in D$.
Define $\tilde{\psi}: M\to \tilde{M}$ by $\tilde{\psi}|_{W^\alpha} =
c(\alpha,0)\psi_\alpha$.
Then, for $a\in V^\alpha$, we have
$$
\begin{array}{l}
  Y^2(a,z)\tilde{\psi}|_{W^\beta} = c(\beta,0) Y^2(a,z) \psi_\beta
  \vsb\\
  = c(\beta,0)c(\alpha,\beta) \psi_{\alpha,\beta}Y^1(a,z)
  \vsb\\
  = c(\alpha +\beta,0) \psi_{\alpha,\beta} Y^1(a,z)\hspace{2cm}
    \text{by}\ \eqref{cocycle}
  \vsb\\
  = \tilde{\psi}|_{W^{\alpha +\beta}} Y^1(a,z).
\end{array}
$$
Therefore, $\tilde{\psi}$ defines a $V_D$-isomorphism between $M$ and
$\tilde{M}$. 
This completes the proof.
\qed

\begin{rem}
  In the case that $D$ is a cyclic group generated by a generator
  $\sigma$, the previous assertion claims that the structure of
  a $\sigma$-stable $V_D$-module is unique over $\C$.
\end{rem}

Next, we consider the fusion rules for simple current extensions.
The following assertion is a direct consequence of the associativity 
\eqref{associativity} for intertwining operators.

\begin{lem}\label{restriction}
  (\cite{DL})
  Let $V_D$ be a $D$-graded extension and let $X$, $W$
  and $T$ be irreducible $V_D$-modules.
  Let $X^0$ and $W^0$ be irreducible $V^0$-submodules of $X$ and $W$, 
  respectively.
  Denote by $\binom{T}{X\ W}_{V_D}$ the space of $V_D$-intertwining
  operators of type $X\times W\to T$.
  Then by a restriction we obtain the following injection:
  $$
    \pi : \binom{T}{X\ W}\ni I(\cd,z)\mapsto I(\cd,z)|_{X^0\tensor
    W^0} \in \binom{T}{X^0\ W^0}_{V^0}.
  $$
\end{lem}

\begin{rem}\label{SCE}
  By the lemma above, we can prove that for any subgroup $E$ of $D$,
  $V_{D/E} = \oplus_{i=1}^{\abs{D/E}} V^{t^i+E}$ is a $D/E$-graded
  simple current extension of $V_E$ if $V_D$ is a $D$-graded simple
  current extension of $V^0$, where $D=\bigcup_{i=1}^{\abs{D/E}}
  (t^i+E)$ denotes a coset decomposition of $D$ with respect to $E$. 
\end{rem}

We prove that the injection $\pi$ becomes an isomorphism in the case
when  $V^0$ contains a tensor product VOA $L(c_{m_1},0)\tensor \cds
\tensor L(c_{m_k},0)$, $V_D$ is a $D$-graded simple current extension
$V^0$ and all of $X$, $W$ and $T$ are $D$-stable. 

\begin{lem}\label{lifting}
  (\cite[Lemma 5.3]{LLY})
  Assume that $V^0$ contains a sub VOA isomorphic to a tensor product
  $L(c_{m_1},0)\tensor \cds \tensor L(c_{m_k},0)$ of unitary Virasoro
  VOAs sharing the same Virasoro vector.
  Assume that $V_D$ is a $D$-graded simple current extension of $V^0$.
  Let $X$, $W$ and $T$ be $D$-stable irreducible $V_D$-modules
  and let $X^0$, $W^0$ and $T^0$ be irreducible $V^0$-submodules 
  of $X$, $W$ and $T$, respectively.
  For any $V^0$-intertwining operator $I(\cd,z)$ of type $X^0\times
  W^0\to T^0$, there exists a $V_D$-intertwining operator
  $\tilde{I}(\cd,z)$ of type $X\times W\to T$ such that
  $\tilde{I}(\cd,z)|_{X^0\tensor W^0} = I(\cd,z)$.
\end{lem}

\pf
The idea of the proof is almost the same as that of Lemma 5.3 of
\cite{LLY}.
In it, they used a result of Huang.
In \cite{H}, Huang proved the following important theorem.

\begin{thm}\label{Huang}
  \cite[Theorem 3.2 and Theorem 3.5]{H}
  Let $V$ be a VOA containing a sub VOA isomorphic to
  $L(c_{m_1},0)\tensor \cds \tensor L(c_{m_k},0)$ with the same
  Virasoro vector.
  Then for any $V$-modules $M^1$, $M^2$, $M^3$, $M^4$, and $M^5$ and
  any intertwining operators $J^1(\cd,z)$ and $J^2(\cd,z)$ of type
  $M^1\times M^4\to M^5$ and $M^2\times M^3\to M^4$, respectively,
  there exist $V$-modules $M^6$ and $M^7$ and intertwining operators
  $J^3(\cd,z)$, $J^4(\cd,z)$, $J^5(\cd,z)$, and $J^6(\cd,z)$ of type
  $M^1\times M^2\to M^6$, $M^6\times M^3\to M^5$, $M^2\times M^7\to
  M^5$, and $M^1\times M^3\to M^7$, respectively, such that
  $\la (w^5)^*, J^1(w^1,z_1)J^2(w^2,z_2)w^3\ra$, $\la (w^5)^*,
  J^4(J^3(w^1,z_1-z_2)w^2,z_2)w^3\ra$ and $\la (w^5)^*,
  J^5(w^2,z_2)J^6(w^1,z_1)w^3\ra$ converge to analytic functions in
  the domains $\abs{z_1}> \abs{z_2}>0$, $\abs{z_2}>\abs{z_1-z_2}>0$ 
  and $\abs{z_2}>\abs{z_1}>0$, respectively, and they are analytic
  extensions of each other after fixing any choices of $\log z_1$,
  $\log z_2$ and $\log (z_1-z_2)$, where $w^i\in M^i$ for $i=1,2,3$
  and $(w^5)^*\in (M^5)^*$. 
\end{thm}

Using the result of Huang above, we prove Lemma \ref{lifting}.
By assumption, we have $D$-graded decompositions $X=\oplus_{\alpha\in
D} X^\alpha$, $W=\oplus_{\alpha\in D} W^\alpha$ and $T=\oplus_{\alpha
\in D}T^\alpha$ such that all $X^\alpha$, $W^\alpha$ and $T^\alpha$,
$\alpha \in D$, are irreducible $V^0$-submodules.
By Theorem \ref{Huang}, there exist $V^0$-intertwining operators
$I^{\alpha,0}(\cd,z)$ and $I^{0,\alpha}(\cd,z)$ of type $X^\alpha
\times W^0 \to T^\alpha$ and $X^0\times W^\alpha\to T^\alpha$,  
respectively such that 
\begin{equation}
  \iota^{-1}_{20}
    \la t^*, I^{\alpha,0}(Y(u^\alpha,z_0)x^0,z_2)w^0
    \ra|_{z_0=z_1-z_2}
  =  \iota^{-1}_{12} \la t^*, Y(u^\alpha,z_1) I^{0,0}(x^0,z_2) 
     w^0 \ra 
\end{equation}
and
\begin{equation}
  \iota^{-1}_{12} \la t^*, Y(u^\alpha,z_1) I^{0,0}(x^0,z_2) w^0 \ra 
  = \iota^{-1}_{21}\la t^*, I^{0,\alpha}(x^0,z_2) Y(u^\alpha,z_1)
    w^0\ra 
\end{equation}
because all $V^\alpha$ are simple current $V^0$-modules,
where $u^\alpha \in V^\alpha$, $x^0\in X^0$, $w^0\in W^0$, $t^*\in
T^*$, and $\iota^{-1}_{12} f(z_1,z_2)$ denotes the formal power
expansion of an analytic function $f(z_1,z_2)$ in the domain
$\abs{z_1} > \abs{z_2}$ (cf. \cite{FHL}).
Then, again by Theorem \ref{Huang}, we can find $V^0$-intertwining
operators $I^{\alpha,\beta}(\cd,z)$ of type $X^\alpha\times W^\beta
\to T^{\alpha+\beta}$ such that
\begin{equation}
  \iota^{-1}_{12} \la t^*, Y(u^\alpha,z_1)I^{0,\beta}(x^0,z_2) 
    w^\beta \ra
  = \iota^{-1}_{20} \la t^*, I^{\alpha,\beta}(Y(u^\alpha,z_0)x^0,z_2) 
    w^\beta \ra |_{z_0=z_1-z_2}.
\end{equation}
We claim that $\tilde{I}(x^\alpha,z) w^\beta :=
I^{\alpha,\beta}(x^\alpha,z) w^\beta$ defines a $V_D$-intertwining
operator of type $X\times W\to T$. 
We only need to show the associativity and the commutativity of
$\tilde{I}(\cd,z)$.
Let $v^\beta\in V^\beta$ and $w^\gamma\in W^\gamma$.
Then we have
$$
\begin{array}{l}
  \iota^{-1}_{120}
    \la t^*, Y(u^\alpha,z_1) I^{\beta,\gamma}(Y(v^\beta,z_0)x^0,z_2) 
    w^\gamma \ra|_{z_0=z_3-z_2}
  \vsb\\
  = \iota^{-1}_{132}
    \la t^*, Y(u^\alpha,z_1) Y(v^\beta,z_3) 
    I^{0,\gamma}(x^0,z_2) w^\gamma\ra
  \vsb\\
  = \iota^{-1}_{342} \la t^*, Y(Y(u^\alpha,z_4)v^\beta,z_3) 
    I^{0,\gamma}(x^0,z_2) w^\gamma\ra|_{z_4=z_1-z_3}
  \vsb\\
  = \iota^{-1}_{240} \la t^*, 
    I^{\alpha+\beta,\gamma}(Y(Y(u^\alpha,z_4)v^\beta,z_0) x^0,z_2) 
    w^\gamma \ra|_{z_4=z_1-z_3,z_0=z_3-z_2}
  \vsb\\
  = \iota^{-1}_{260} \la t^*, 
    I^{\alpha +\beta,\gamma}(Y(u^\alpha,z_6) Y(v^\beta,z_0)x^0,z_2) 
    w^\gamma\ra|_{z_6=z_1-z_2,z_0=z_3-z_2}
\end{array}
$$
and hence we obtain the following associativity:
\begin{equation} 
  \la t^*, Y(u^\alpha,z_1)I^{\beta,\gamma}(x^\beta,z_2) w^\gamma \ra
    = \la t^*, I^{\alpha+\beta,\gamma}(Y(u^\alpha,z_0)x^\beta,z_2)
      w^\gamma\ra|_{z_0=z_1-z_2}.
\end{equation}

Next we prove the commutativity of $I^{\alpha,\beta}(\cd,z)$.
We have
$$
\begin{array}{l}
  \iota^{-1}_{201} \la t^*, I^{\beta,\alpha}(Y(v^\beta,z_0)x^0,z_2)
    Y(u^\alpha,z_1)w^0\ra|_{z_0=z_3-z_2}
  \vsb\\
  = \iota^{-1}_{321} \la t^*, Y(v^\beta,z_3) I^{0,\alpha}(x^0,z_2)
    Y(u^\alpha,z_1)w^0\ra
  \vsb\\
  = \iota^{-1}_{312} \la t^*, Y(v^\beta,z_3) Y(u^\alpha,z_1) 
    I^{0,0}(x^0,z_2) w^0\ra
\end{array}
$$

$$
\begin{array}{l}
  = \iota^{-1}_{132} \la t^*, Y(u^\alpha,z_1) Y(v^\beta,z_3)
    I^{0,0}(x^0,z_2) w^0\ra
  \vsb\\
  = \iota^{-1}_{342} \la t^*, Y(Y(u^\alpha,z_4)v^\beta,z_3)
    I^{0,0}(x^0,z_2) w^0\ra|_{z_4=z_1-z_3}
  \vsb\\
  = \iota^{-1}_{204} \la t^*, I^{\alpha+\beta,0}(Y(Y(u^\alpha,z_4)
    v^\beta,z_0)x^0,z_2) w^0\ra|_{z_0=z_3-z_2,z_4=z_1-z_3}
  \vsb\\
  = \iota^{-1}_{250} \la t^*, I^{\alpha+\beta,0}(Y(u^\alpha,z_5)
    Y(v^\beta,z_0)x^0,z_2) w^0\ra|_{z_0=z_3-z_2, z_5=z_1-z_2}
  \vsb\\
  = \iota_{120} \la t^*, Y(u^\alpha,z_1) 
    I^{\beta,0}(Y(v^\beta,z_0)x^0,z_2) w^0\ra|_{z_0=z_3-z_2}.
\end{array}
$$
Thus, we get the following:
\begin{equation}
  \la t^*, Y(u^\alpha,z_1) I^{\beta,0}(x^\beta,z_2) w^0\ra
  = \la t^*, I^{\beta,\alpha}(x^\beta,z_2) Y(u^\alpha,z_1) w^0\ra .
\end{equation}
Then
$$
\begin{array}{l}
  \iota^{-1}_{123} \la t^*, Y(u^\alpha,z_1)
    I^{\beta,\gamma}(x^\beta,z_2) Y(v^\gamma,z_3) w^0 \ra 
  \vsb\\
  = \iota^{-1}_{132} \la t^*, Y(u^\alpha,z_1) Y(v^\gamma,z_3)
    I^{\beta,0}(x^\beta,z_2) w^0 \ra
  \vsb\\
  = \iota^{-1}_{302} \la t^*, Y(Y(u^\alpha,z_0)v^\gamma,z_3)
    I^{\beta,0}(x^\beta,z_2) w^0 \ra|_{z_0=z_1-z_3}
  \vsb\\
  = \iota^{-1}_{230} \la t^*, I^{\beta,\alpha+\gamma}(x^\beta,z_2) 
    Y(Y(u^\alpha,z_0)v^\gamma,z_3) w^0 \ra|_{z_0=z_1-z_3}
  \vsb\\
  = \iota^{-1}_{213} \la t^*, I^{\beta,\alpha+\gamma}(x^\beta,z_2)
    Y(u^\alpha,z_1) Y(v^\beta,z_3) w^0 \ra
\end{array}
$$
and hence we arrive at the following commutativity:
\begin{equation}
  \la t^*, Y(u^\alpha,z_1) I^{\beta,\gamma}(x^\beta,z_2) w^\gamma \ra
  = \la t^*, I^{\beta,\alpha +\gamma}(x^\beta,z_2)
    Y(u^\alpha,z_1) w^\gamma \ra .
\end{equation}
This completes the proof of Lemma \ref{lifting}.
\qed

In the rest of this section, we study a relation between automorphisms 
of $V^0$ and those of $V_D$.
Let $\sigma$ be an automorphism of $V^0$ and denote by
$(V^\alpha)^\sigma$ the $\sigma$-conjugate $V^0$-module of $V^\alpha$
for $\alpha \in D$.
If there exists a $D$-graded extension $V_D=\oplus_{\alpha\in D}
V^\alpha$ of $V^0$, then we can construct another $D$-graded extension 
$V_D' = \oplus_{\alpha\in D} (V^\alpha)^\sigma$ in the following way.
By definition, there exist linear isomorphisms $\varphi_\alpha :
V^\alpha \to (V^\alpha)^\sigma$ such that $Y_{(V^\alpha)^\sigma} 
(a,z) \varphi_\alpha = \varphi_\alpha Y_{V^\alpha} (\sigma a,z)$
for all $a\in V^0$. 
For $a\in V^\alpha$ and $b\in V^\beta$, define the vertex operation in 
$V_D'=\oplus_{\alpha\in D}(V^\alpha)^\sigma$ by
$$
  Y_{V_D'}(\varphi_\alpha a,z)\varphi_\beta b
  := \varphi_{\alpha+\beta}Y_{V_D}(a,z)b.
$$
Since $Y_{V_D'} (\cd ,z)|_{(V^\alpha)^\sigma \times (V^\beta)^\sigma}$
is a $V^0$-intertwining operator of type $(V^\alpha)^\sigma \times
(V^\beta)^\sigma \to (V^{\alpha+\beta})^\sigma$,
$(V_D',Y_{V_D'}(\cd,z))$ also forms a $D$-graded extension of $V^0$.
Moreover, if $V_D$ is a $D$-graded simple current extension of $V^0$,
then so is $V_D'$.
We call $V_D'$ the {\it $\sigma$-conjugate} of $V_D$.
It is clear from its construction that $V_D$ and $V_D'$ are
isomorphic as VOAs even if $\{ V^\alpha \mid \alpha\in D\}$ and 
$\{ (V^\alpha)^\sigma \mid \alpha\in D\}$ are distinct sets of
inequivalent $V^0$-modules.
Therefore, we introduce the following definition.

\begin{df}
  Two $D$-graded simple current extensions $V_D=\oplus_{\alpha\in
  D}V^\alpha$ and $\tilde{V}_D = \oplus_{\alpha\in
  D}\tilde{V}^\alpha$ are said to be {\it equivalent} if there
  exists a VOA-isomorphism $\Phi : V_D\to \tilde{V}_D$ such that
  $\Phi (V^\alpha) = \tilde{V}^\alpha$ for all $\alpha \in D$.
\end{df}

The following assertion will be needed later.

\begin{lem}\label{stability}
  Suppose that $V_D$ is a $D$-graded extension of $V^0$. 
  For an automorphism $\sigma \in \aut (V^0)$, assume that there is an
  automorphism $\Psi$ on $V_D$ such that $\Psi 
  (V^0)=V^0$ and $\Psi|_{V^0}=\sigma$.
  Then as sets of inequivalent irreducible $V^0$-modules, 
  $\{ \Psi^{-1} V^\alpha \mid \alpha \in D\}$ and 
  $\{ (V^\alpha)^\sigma \mid \alpha \in D\}$ are the same.
\end{lem}

\pf
Denote $Y_{V_D}(\cd,z)|_{V^0\tensor V^\alpha}$ by $Y_\alpha (\cd,z)$.
By definition, we can take linear isomorphisms $\varphi_\alpha :
V^\alpha \to (V^\alpha)^\sigma$ such that $Y_{(V^\alpha)^\sigma}(a,z)
\varphi_\alpha = \varphi_\alpha Y_\alpha (\sigma a,z)$ for all $a\in
V^0$. 
Define $\Psi_\alpha : \Psi^{-1} V^\alpha \to (V^\alpha)^\sigma$ by
$\Psi_\alpha = \varphi_\alpha \circ \Psi|_{\Psi^{-1}V^\alpha}$.
Then for $a\in V^0$ we have 
$$
\begin{array}{l}
  Y_{(V^\alpha)^\sigma}(a,z) \Psi_\alpha 
  = Y_{(V^\alpha)^\sigma}(a,z) \varphi_\alpha \Psi =\varphi_\alpha 
    Y_\alpha (\sigma a,z)\Psi 
  \vsb\\
  = \varphi_\alpha Y_\alpha (\Psi a,z) \Psi 
  = \varphi_\alpha \Psi Y_{V_D}(a,z)|_{\Psi^{-1}V^\alpha}
  = \Psi_\alpha Y_{V_D}(a,z)|_{\Psi^{-1} V^\alpha}.
\end{array}
$$
Therefore, $\Psi_\alpha$ is a $V^0$-isomorphisms.
Hence, we get the assertion.
\qed

\section{Vertex operator algebra with two Miyamoto involutions
generating $S_3$} 

In this section we study a VOA on which $S_3$ acts.
First, we construct it from a lattice VOA.
More precisely, we will find it in an extension of an affine VOA.
Then we show that there exists a unique VOA structure on it.
All irreducible modules are classified.
At last, we prove that they are generated by two conformal vectors
with central charge $1/2$ and the full automorphism group is
isomorphic to $S_3$.
Namely, it is the VOA of involution type $A_2$ in the sense of
Miyamoto \cite{M3}.

\subsection{Construction} 

Let $A_1^5=\Z \alpha^1 \oplus \Z \alpha^2 \oplus \cds \oplus \Z
\alpha^5$ with $\la \alpha^i,\alpha^j\ra =2\delta_{i,j}$ and 
set $L:=A_1^5\cup (\gamma +A_1^5)$ with $\gamma := \hf \alpha^1 +\hf
\alpha^2 +\hf \alpha^3 +\hf \alpha^4$. 
Then $L$ is an even lattice so that we can construct a VOA $V_L$
associated to $L$.
We have an isomorphism $V_L= V_{A_1^5}\oplus V_{\gamma +A_1^5}\simeq
\{ \mL (1,0)^{\tensor 4}\oplus \mL (1,1)^{\tensor 4}\} \tensor \mL
(1,0)$. 
By \eqref{character} and the fusion rules \eqref{wang} and \eqref{FZ}, 
we can show the following.

\begin{lem}
  We have the following inclusions
  $$
  \begin{array}{lll}
  \mL (1,0)^{\tensor 3} 
    &\supset&  L(\fr{1}{2},0) \tensor L(\fr{7}{10},0) \tensor \mL (3,0),
  \vsb\\
  \mL (1,1)^{\tensor 3} &\supset& L(\fr{1}{2},0) \tensor
    L(\fr{7}{10},0) \tensor \mL (3,3) .
  \end{array}
  $$
  Therefore, $V_L$ contains a sub VOA isomorphic to
  $$
    \mL (3,0)\tensor \mL (1,0)\tensor \mL (1,0)
    \oplus 
    \mL(3,3)\tensor \mL (1,1) \tensor \mL (1,0) .
  $$ 
\end{lem}

\begin{lem}\label{existence}
  We have the following decompositions:
  $$
  \begin{array}{l}
  \begin{array}{lll}
    \begin{array}{l}
    \mL (3,0) \tensor \mL (1,0) \tensor \mL (1,0)
    \vsb\\
    \simeq \l\{
      \begin{array}{c}
        L(\fr{4}{5},0)\tensor L(\fr{6}{7},0)
        \\ \oplus \\
        L(\fr{4}{5},3)\tensor L(\fr{6}{7},5)
        \\ \oplus\\
        L(\fr{4}{5},\fr{2}{3})\tensor L(\fr{6}{7},\fr{4}{3})
      \end{array}
    \r\} \tensor \mL (5,0) 
    \vsb\\
    \bigoplus
    \l\{
      \begin{array}{c}
        L(\fr{4}{5},0)\tensor L(\fr{6}{7},\fr{5}{7})
        \\ \oplus \\
        L(\fr{4}{5},3)\tensor L(\fr{6}{7},\fr{12}{7})
        \\ \oplus \\
        L(\fr{4}{5},\fr{2}{3})\tensor L(\fr{6}{7},\fr{1}{21})
      \end{array}
    \r\} \tensor \mL (5,2) 
    \vsb\\
    \bigoplus
    \l\{
      \begin{array}{c}
        L(\fr{4}{5},0)\tensor L(\fr{6}{7},\fr{22}{7})
        \\ \oplus \\
        L(\fr{4}{5},3)\tensor L(\fr{6}{7},\fr{1}{7})
        \\ \oplus \\
        L(\fr{4}{5},\fr{2}{3})\tensor L(\fr{6}{7},\fr{10}{21})
      \end{array}
    \r\} \tensor \mL (5,4),
    \end{array}
    &
    \q
    &
    \begin{array}{l}
    \mL (3,3) \tensor \mL (1,1)\tensor \mL (1,0)
    \vsb\\
    \simeq \l\{
      \begin{array}{c}
        L(\fr{4}{5},0)\tensor L(\fr{6}{7},5)
        \\ \oplus \\
        L(\fr{4}{5},3)\tensor L(\fr{6}{7},0)
        \\ \oplus \\
        L(\fr{4}{5},\fr{2}{3})\tensor L(\fr{6}{7},\fr{4}{3})
      \end{array}
    \r\} \tensor \mL (5,0)
    \vsb\\
    \bigoplus \l\{
      \begin{array}{c}
        L(\fr{4}{5},0)\tensor L(\fr{6}{7},\fr{12}{7})
        \\ \oplus \\
        L(\fr{4}{5},3)\tensor L(\fr{6}{7},\fr{5}{7})
        \\ \oplus \\
        L(\fr{4}{5},\fr{2}{3})\tensor L(\fr{6}{7},\fr{1}{21})
      \end{array}
      \r\} \tensor \mL (5,2)
    \vsb\\
    \bigoplus \l\{
      \begin{array}{c}
        L(\fr{4}{5},0)\tensor L(\fr{6}{7},\fr{1}{7})
        \\ \oplus \\
        L(\fr{4}{5},3)\tensor L(\fr{6}{7},\fr{22}{7})
        \\ \oplus \\
        L(\fr{4}{5},\fr{2}{3})\tensor L(\fr{6}{7},\fr{10}{21})
     \end{array}
     \r\} \tensor \mL (5,4).
  \end{array}
  \end{array}
  \end{array}
  $$
  Hence, $\mL (3,0)\tensor \mL (1,0)\tensor \mL (1,0) \oplus \mL
  (3,3)\tensor \mL (1,1) \tensor \mL (1,0)$ (and $V_L$) contains a sub
  VOA $U$ isomorphic to 
  \begin{equation}\label{U}
  \begin{array}{lll}
    \begin{bmatrix}
      \begin{array}{c}
        L(\fr{4}{5},0)\tensor L(\fr{6}{7},0)
        \\ \oplus \\
        L(\fr{4}{5},3)\tensor L(\fr{6}{7},5)
        \\ \oplus \\
        L(\fr{4}{5},\fr{2}{3})\tensor L(\fr{6}{7},\fr{4}{3})
      \end{array}
    \end{bmatrix}
    & 
    \bigoplus
    &
    \begin{bmatrix}
      \begin{array}{c}
        L(\fr{4}{5},0)\tensor L(\fr{6}{7},5)
        \\ \oplus \\
        L(\fr{4}{5},3)\tensor L(\fr{6}{7},0)
        \\ \oplus \\
        L(\fr{4}{5},\fr{2}{3})\tensor L(\fr{6}{7},\fr{4}{3})
      \end{array}
    \end{bmatrix}
  \end{array}.
  \end{equation}
\end{lem}

\begin{rem}
  Note that the sub VOA $U$ has exactly the same form as stated in 
  Theorem 5.6 (4) of \cite{M3}.
  In the following context, we will show that our VOA $U$ is actually 
  the same as $\text{VA}(e,f)$ in \cite{M3}.
\end{rem}

\begin{rem}
  By the lemma above, we note that $L(\fr{6}{7},0)\oplus
  L(\fr{6}{7},5)$ is a sub VOA of $U$, which completes the proof of
  Theorem \ref{tricritical}. 
\end{rem}

We can also define $U$ in the following way.
For $i=1,2,\dots,5$, set
$$
\begin{array}{l}
  H^j:= \alpha^1_{(-1)}\vac + \cds + \alpha^j_{(-1)}\vac,
  \vsb\\
  E^j:= e^{\alpha^1} + \cds + e^{\alpha^j},
  \vsb\\
  F^j:= e^{-\alpha^1} + \cds + e^{-\alpha^j},
  \vsb\\
  \Omega^j := \dfr{1}{2(j+2)}\l( \dfr{1}{2} H^j_{(-1)} H^j 
    + E^j_{(-1)} F^j + F^j_{(-1)} E^j \r) ,\vsb\\
  \w^i := \Omega^i +\dfr{1}{4}\l(\alpha^{i+1}_{(-1)}\r)^2 \vac 
    - \Omega^{i+1}.
\end{array}
$$
Then $H^j,E^j$ and $F^j$ generate a simple affine sub VOA $\mL (j,0)$
and $\w^i$, $1\leq i\leq 4$, generate simple Virasoro sub VOAs  
$L(c_i,0)$ in $V_L$.
Furthermore, we have an orthogonal decomposition of the
Virasoro vector $\w_{V_L}$ of $V_L$ into a sum of mutually commutative
Virasoro vectors as
$$
  \w_{V_L} = \w^1 + \w^2 + \w^3 + \w^4 +\Omega^5.
$$
Then we may define $U$ to be as follows:
$$
  U = \{ v\in V_L \mid \w^1_{(1)} v = \w^2_{(1)} v = \Omega^5_{(1)} 
    v =0\} .
$$
Set
\begin{equation}\label{v_0}
\begin{array}{ll}
  e &:= \dfr{1}{16}\l((\alpha^4 -\alpha^5)_{(-1)}\r)^2\vac 
  -\dfr{1}{4}(e^{\alpha^4-\alpha^5}+ e^{-\alpha^4 +\alpha^5}),
  \vsb\\
  v_0 &:= \dfr{5}{18}\w^3+\dfr{7}{9}\w^4 -\dfr{16}{9}e,
  \vsb\\
  v_1 &:= (9F^4- 8F^5)_{(-1)} (4F^3-3F^4)_{(0)} 
    e^{\fr{1}{2}(\alpha^1 + \alpha^2  +\alpha^3 +\alpha^4)}
  \vsb\\
    &\q\qq -\dfr{1}{2} (9H^4-8H^5)_{(-1)} F^4_{(0)} 
      (4F^3-3F^4)_{(0)} 
      e^{\fr{1}{2}(\alpha^1 + \alpha^2  +\alpha^3 +\alpha^4)}
  \vsb\\
    & \q\qq\qq -\dfr{1}{2} (9E^4-8E^5)_{(-1)} \l( F^4_{(0)}\r)^2 
      (4F^3-3F^4)_{(0)} e^{\fr{1}{2}(\alpha^1 + \alpha^2  +\alpha^3
      +\alpha^4)}. 
\end{array}
\end{equation}

Then we can show that both $e$ and $v_i$ are contained in $U_2$ and
$e_{(1)}e=2e$, $e_{(3)}e=\fr{1}{4}\vac$, $\w^3_{(1)} v_i=\fr{2}{3}v_i$ 
and $\w^4_{(1)} v_i=\fr{4}{3}v_i$ for $i=0,1$. 
Therefore, $e$ generates a sub VOA isomorphic to $L(\fr{1}{2},0)$ in
$U$ and $v_i$, $i=0,1$, are highest weight vectors for $\la \w^3\ra
\tensor \la \w^4\ra \simeq L(\fr{4}{5},0)\tensor L(\fr{6}{7},0)$ with
highest weight $(\fr{2}{3},\fr{4}{3})$.
Since the weight 2 subspace of $U$ is 4-dimensional, we note
that $\w^3,\w^4, v_0$ and $v_1$ span $U_2$.
In the next subsection we will show that they generate $U$ as a
VOA.

\subsection{Structures} 

By Lemma \ref{existence}, we know that there exists a structure of a
VOA in \eqref{U}.
Here we will prove that there exists a unique VOA structure on it.
By \eqref{U}, $U$ contains a tensor product of two extensions of the
unitary Virasoro VOAs $W(0)=L(\fr{4}{5},0)\oplus L(\fr{4}{5},3)$ and
$N(0)=L(\fr{6}{7},0)\oplus L(\fr{6}{7},5)$ (See Sec.{} 2.4--5). 
Since both $W(0)$ and $N(0)$ are rational, $U$ is completely reducible 
as a $W(0)\tensor N(0)$-module.
Therefore, $U$ as a $W(0)\tensor N(0)$-module is isomorphic to
$$
  U \simeq W(0)\tensor N(0) 
  \oplus W\l(\fr{2}{3}\r)^{\epsilon_1}\tensor N\l(\fr{4}{3}\r)^{\xi_1} 
  \oplus W\l(\fr{2}{3}\r)^{\epsilon_2}\tensor
    N\l(\fr{4}{3}\r)^{\xi_2}, 
$$
where $\epsilon_i,\xi_j=\pm$.
Recall that both $W(0)$ and $N(0)$ have the canonical involutions
$\sigma_1$ and $\sigma_2$, respectively.
Then they can be lifted to involutions of $W(0)\tensor N(0)$
and we still denote them by $\sigma_1$ and $\sigma_2$, respectively.
By our construction, $U$ has a $\Z_2$-grading $U=U^+\oplus U^-$ with 
$U^+\subset \mL (3,0)\tensor \mL (1,0)\tensor \mL (1,0)\subset
V_{A_1^5}$ and $U^-\subset \mL (3,3)\tensor \mL (1,1) \tensor \mL
(1,0) \subset V_{\gamma +A_1^5}$.
We note that this decomposition defines a natural extension of
involution $\sigma_1 \sigma_2$ on $W(0)\tensor N(0)$ to that on
$U$, which we will also denote by $\sigma_1 \sigma_2$.
Therefore, by Lemma \ref{stability} we have
$(W(\fr{2}{3})^{\epsilon_1}\tensor 
N(\fr{4}{3})^{\xi_1})^{\sigma_1 \sigma_2} =
W(\fr{2}{3})^{\epsilon_2} \tensor N(\fr{4}{3})^{\xi_2}$
and hence $\epsilon_2=-\epsilon_1$ and $\xi_2=-\xi_1$.
Since we may rename the signs of the irreducible $N(0)$-modules of
$\pm$-type (cf. Remark \ref{signature}), we may assume that
$\epsilon_1=\xi_1$. 

\begin{thm}\label{gen}
  A VOA $U$ contains a sub VOA $W(0)\tensor N(0)$. 
  As a $W(0)\tensor N(0)$-module, $U$ is isomorphic to 
  \begin{equation}\label{(++)-type}
    W(0)\tensor N(0)
    \oplus W\l(\fr{2}{3}\r)^+ \tensor N\l(\fr{4}{3}\r)^+
    \oplus W\l(\fr{2}{3}\r)^- \tensor N\l(\fr{4}{3}\r)^-
  \end{equation}
  after fixing suitable choice of $\pm$-type of $N(\fr{4}{3})^\pm$. 
  Therefore, $U$ is a simple VOA and generated by its weight 2 
  subspace as a VOA.
\end{thm}

\pf
The decomposition is already shown.
Since $U$ is a $\Z_3$-simple current extension of $W(0)\tensor N(0)$, 
it is a simple VOA.
So we should show that $U_2$ generates $U$.
Since $U_2$ contains the Virasoro vectors $\w^3$ and $\w^4$ and
highest weight vectors of $W(\fr{2}{3})^\pm \tensor N(\fr{4}{3})^\pm$,
$U_2$ generates whole of $W(\fr{2}{3})^\pm\tensor N(\fr{4}{3})^\pm$.
Since $V_L$ is simple, for any non-zero vectors $u\in
W(\fr{2}{3})^+\tensor N(\fr{4}{3})^+$ and $v\in W(\fr{2}{3})^-\tensor
N(\fr{4}{3})^-$ we have $Y(u,z)v\ne 0$ in $U$ (cf. \cite{DL}).
Therefore, by the fusion rules in Theorem \ref{fusion for W(0)} and 
Theorem \ref{fusion for N(0)}, $W(\fr{2}{3})^\pm \tensor
N(\fr{4}{3})^\pm$ generate $W(0)\tensor N(0)$ in $U$.
Hence, $U_2$ generates whole of $U$.
\qed

By the previous theorem together with the arguments in Section 3, we
note that there exists the following $\Z_3$-simple current extension
of $W(0)\tensor N(0)$.
\begin{equation}\label{(+-)-type}
  U'=W(0)\tensor N(0)
  \oplus W\l(\fr{2}{3}\r)^+\tensor N\l(\fr{4}{3}\r)^-
  \oplus W\l(\fr{2}{3}\r)^-\tensor N\l(\fr{4}{3}\r)^+.
\end{equation}

Since $U$ and $U'$ are $\sigma_1$-conjugate extensions of each others,
they are equivalent $\Z_3$-simple current extensions of $W(0)\tensor
N(0)$. 
Thus, we get the following.

\begin{thm}
  The following $\Z_3$-simple current extensions of $W(0)\tensor N(0)$ 
  are equivalent:
  $$
    W(0)\tensor N(0)\oplus W\l(\fr{2}{3}\r)^+ \tensor
    N\l(\fr{4}{3}\r)^\pm \oplus W\l(\fr{2}{3}\r)^- \tensor 
    N\l(\fr{4}{3}\r)^\mp .
  $$
  Hence, there is a unique $\Z_3$-graded VOA structure in \eqref{U}. 
\end{thm}

\subsection{Modules} 

Let $U$ be the $\Z_3$-graded VOA as in \eqref{U}.
In this subsection we will classify all irreducible $U$-modules.
Set $U=U^0\oplus U^1\oplus U^2$ with $U^0=W(0)\tensor N(0)$,
$U^1=W(\fr{2}{3})^+\tensor N(\fr{4}{3})^+$ and
$U^2=W(\fr{2}{3})^-\tensor N(\fr{4}{3})^-$. 

\begin{lem}
  Every irreducible $U$-modules is $\Z_3$-stable.
\end{lem}

\pf
Let $M$ be an irreducible $U$-module.
Take an irreducible $U^0$-submodule $P$ of $M$.
By Lemma \ref{non-zero}, both $U^1\cd P$ and $U^2\cd P$ are non-zero
irreducible $U^0$-submodules of $M$.
It follows from the fusion rules for $U^0=W(0)\tensor N(0)$-modules
that $U^i\cd P \not\simeq U^j\cd P$ as $U^0$-modules if $i\not\equiv j 
\mod 3$.
Therefore, $M=P\oplus (U^1\cd P)\oplus (U^2\cd P)$ and hence 
$M$ has a $\Z_3$-grading under the action of $U$.
This completes the proof.
\qed

By this lemma and Proposition \ref{unique1}, the $U^0$-module 
structure of each irreducible $U$-module completely determines the 
$U$-module structure of it. 

\begin{lem}\label{candidates}
  Let $M$ be an irreducible $U$-module. 
  Then, as a $W(0)\tensor N(0)$-module, $M$ is isomorphic to one of
  the following:  
  $$
  \begin{array}{l}
    W(0)\tensor N(0) 
      \oplus W(\fr{2}{3})^+\tensor N(\fr{4}{3})^+ 
      \oplus W(\fr{2}{3})^-\tensor N(\fr{4}{3})^-,
    \vsb\\
    W(0)\tensor N(\fr{1}{7})
      \oplus W(\fr{2}{3})^+\tensor N(\fr{10}{21})^+ 
      \oplus W(\fr{2}{3})^- \tensor N(\fr{10}{21})^-,
    \vsb\\
    W(0)\tensor N(\fr{5}{7})
      \oplus W(\fr{2}{3})^+\tensor N(\fr{1}{21})^+
      \oplus W(\fr{2}{3})^-\tensor N(\fr{1}{21})^-,
    \vsb\\
    W(\fr{2}{5})\tensor N(0)
      \oplus W(\fr{1}{15})^+\tensor N(\fr{4}{3})^+
      \oplus W(\fr{1}{15})^-\tensor N(\fr{4}{3})^-,
    \vsb\\  
    W(\fr{2}{5})\tensor N(\fr{1}{7})
      \oplus W(\fr{1}{15})^+\tensor N(\fr{10}{21})^+
      \oplus W(\fr{1}{15})^-\tensor N(\fr{10}{21})^-,
    \vsb\\
    W(\fr{2}{5})\tensor N(\fr{5}{7})
      \oplus W(\fr{1}{15})^+\tensor N(\fr{1}{21})^+
      \oplus W(\fr{1}{15})^-\tensor N(\fr{1}{21})^-.
  \end{array}
  $$
\end{lem}

\pf
Let $M$ be an irreducible $U$-module and $P^0$ an irreducible
$U^0$-submodule of $M$.
Then $M=P^0\oplus P^1 \oplus P^2$ with $P^1=U^1\times P^0$ and
$P^2=U^2\times P^0$.
The vertex operators $Y_M(\cd,z)$ on $M$ give $U^0$-intertwining
operators of type $U^i\times P^j\to P^{i+j}$ for $i,j\in\Z_3$.
The powers of $z$ in an intertwining operator of type $U^i\times
P^j\to P^{i+j}$ are contained in $-h_{U^i}-h_{P^i}+h_{P^{i+j}}+\Z$,
where $h_X$ denotes the top weight of a $U^0$-module $X$.
Since the powers of $z$ in $Y_M(\cd ,z)$ belong to $\Z$, by
considering top weights we arrive at the candidates above.
\qed

\begin{thm}\label{irrs}
  All irreducible $U$-modules are given by the listed in Lemma
  \ref{candidates}.
  That means, they are exactly $U$-modules.
\end{thm}

\pf
We already know that if there exist $U$-module structures in the
candidates in Lemma \ref{candidates}, then they must be unique by
Proposition \ref{unique1}. 
So we only need to show that they are actually $U$-modules.
Recall that $U\tensor \mL (5,0)$ is a sub VOA of a VOA 
$$
  T= \mL (3,0)\tensor \mL (1,0)\tensor \mL (1,0) 
  \oplus \mL (3,3)\tensor \mL (1,1)\tensor \mL (1,0).
$$
It is shown in \cite{Li2} that 
$$
  \mL (3,2)\tensor \mL (1,0) \oplus 
  \mL (3,1) \tensor \mL (1,1)
$$ is an irreducible  $\mL (3,0)\tensor \mL
(1,0)$ $\oplus \mL (3,3)\tensor \mL (1,1)$-module.
Hence, 
$$
  \mL (3,2)\tensor \mL (1,0)\tensor \mL (1,0)
  \oplus 
  \mL (3,1)\tensor \mL (1,1)\tensor \mL (1,0)
$$ is an irreducible $T$-module.
Then by using \eqref{character}, we get the following decompositions:
$$
\begin{array}{l}
  \mL (3,0)\tensor \mL (1,0)\tensor \mL (1,0) 
    \oplus \mL (3,3)\tensor \mL (1,1)\tensor \mL (1,0)
    \vsb\\
    \simeq
    \{ W(0)\tensor N(0) 
      \oplus W(\fr{2}{3})^+\tensor N(\fr{4}{3})^+ 
      \oplus W(\fr{2}{3})^-\tensor N(\fr{4}{3})^-
    \}
    \tensor \mL (5,0)
    \vsb\\
    \q \oplus
    \{ W(0)\tensor N(\fr{5}{7})
      \oplus W(\fr{2}{3})^+\tensor N(\fr{1}{21})^+
      \oplus W(\fr{2}{3})^-\tensor N(\fr{1}{21})^-
    \}
    \tensor \mL (5,2)
    \vsb\\
    \q \oplus
    \{ W(0)\tensor N(\fr{1}{7})
      \oplus W(\fr{2}{3})^+\tensor N(\fr{10}{21})^+ 
      \oplus W(\fr{2}{3})^- \tensor N(\fr{10}{21})^-
    \}
    \tensor \mL (5,4),
  \vsv\\
  \mL (3,2)\tensor \mL (1,0)\tensor \mL (1,0)
    \oplus \mL (3,1)\tensor \mL (1,1)\tensor \mL (1,0)
    \vsb\\
    \simeq
    \{ W(\fr{2}{5})\tensor N(0)
      \oplus W(\fr{1}{15})^+\tensor N(\fr{4}{3})^+
      \oplus W(\fr{1}{15})^-\tensor N(\fr{4}{3})^-
    \}
    \tensor \mL (5,0)
    \vsb\\
    \q \oplus
    \{ W(\fr{2}{5})\tensor N(\fr{5}{7})
      \oplus W(\fr{1}{15})^+\tensor N(\fr{1}{21})^+
      \oplus W(\fr{1}{15})^-\tensor N(\fr{1}{21})^-
    \}
    \tensor \mL (5,2)
    \vsb\\  
    \q \oplus
    \{ W(\fr{2}{5})\tensor N(\fr{1}{7})
      \oplus W(\fr{1}{15})^+\tensor N(\fr{10}{21})^+
      \oplus W(\fr{1}{15})^-\tensor N(\fr{10}{21})^-
    \}
    \tensor \mL (5,4).
\end{array}
$$
Therefore, all candidates in Lemma \ref{candidates} are $U$-modules.
\qed

\begin{thm}
  $U$ is rational.
\end{thm}

\pf
Let $M$ be an admissible $U$-module.
Take an irreducible $U^0$-submodule $P$.
By Lemma \ref{non-zero}, both $U^1\cd P$ and $U^2\cd P$ are 
non-trivial irreducible $U^0$-submodule of $M$. 
Since $U^i\cd P\not\simeq U^j\cd P$ if $i\not\equiv j\mod 3$, 
$P+(U^1\cd P)+(U^2\cd P)=P\oplus (U^1\cd P)\oplus (U^2\cd P)$ is an
irreducible $U$-submodule of $M$.
Hence, every irreducible $U^0$-submodule of $M$ is contained in an
irreducible $U$-submodule.
Thus $M$ is a completely reducible $U$-module.
\qed

\subsection{Conformal vectors} 

In this subsection we study the Griess algebra of $U$.
Recall $e$, $v_0$, $v_1 \in U_2$ defined  by \eqref{v_0}.
Set
$$
\begin{array}{l}
  \w := \w^3+\w^4,
  \vsb\\
  a := \dfr{105}{2^8} (\w -e),
  \vsb\\
  b := \dfr{3^2}{2^8}(-5\w^3+7\w^4 -4e),
  \vsb\\
  c := k v_1,
\end{array}
$$
where the scalar $k\in \R$ is determined by the condition 
$\la c,c\ra =3^5/2^{11}$. 
Then $\{ e, a, b, c\}$ is a set of basis of $U_2$.
By direct calculations one can show that the multiplications and inner 
products in the Griess algebra of $U$ are given as follows:
$$
\begin{array}{lll}
  e_{(1)} a= 0,
  & e_{(1)} b= \dfr{1}{2}b,
  & e_{(1)} c= \dfr{1}{16}c,
  \vsb\\
  a_{(1)} a= \dfr{105}{2^7}a, 
  & a_{(1)} b= \dfr{3^2\cd 5\cd 7}{2^9}b, 
  & a_{(1)} c= \dfr{31\cd 105}{2^{12}}c,
  \vsb\\
  b_{(1)} b= \dfr{3^7}{2^{15}}e+\dfr{3^3}{2^7}a, 
  & b_{(1)} c= \dfr{3^2\cd 23}{2^{10}} c,
  & c_{(1)} c= \dfr{3^5}{2^{13}}e+\dfr{31}{2^5}a+\dfr{23}{2^5}b,
  \vsb\\
  \la a,a\ra = \dfr{3^6\cd 5\cd 7}{2^{18}},
    & \la b,b\ra = \dfr{3^7}{2^{16}},
    & \la c,c\ra = \dfr{3^5}{2^{11}}.
\end{array}
$$
Hence, we note that the Griess algebra of our VOA $U$ is isomorphic to
that of $\text{VA}(e,f)$ with $\la e,f\ra =13/2^{10}$ in \cite{M3}.
Therefore, by tracing calculations in \cite{M3} we can find 
the following conformal vectors with central charge $1/2$ in $U_2$.
$$
\begin{array}{cc}
  f:= \dfr{13}{2^8}e + a + b + c,
  &
  f':= \dfr{13}{2^8}e + a + b -c.
\end{array}
$$
And by a calculation we get
\begin{equation}\label{coc}
\begin{array}{l}
  e_{(1)} f=-\dfr{105}{2^9} \w +\dfr{9}{2^5}e + \dfr{9}{2^5}f
    + \dfr{7}{2^5}f',
  \vsb\\
  e_{(1)} f' =-\dfr{105}{2^9} \w +\dfr{9}{2^5}e + \dfr{7}{2^5}f
    +\dfr{9}{2^5}f',
  \vsb\\
  f_{(1)} f'=-\dfr{105}{2^9} \w +\dfr{7}{2^5}e + \dfr{9}{2^5}f
    +\dfr{9}{2^5}f',
  \vsb\\
  \la e,f\ra =\la e,f'\ra = \la f,f'\ra =\dfr{13}{2^{10}}.
\end{array}
\end{equation}
Using these equalities, we can show that the Griess algebra $U_2$
is generated by two conformal vectors $e$ and $f$.
Since $U_2$ generates $U$ as a VOA, $U$ is generated by two
conformal vectors $e$ and $f$. Thus

\begin{thm}
  $U$ is generated by two conformal vectors $e$ and $f$ with
  central charge $1/2$ such that $\la e,f\ra =13/2^{10}$.
\end{thm}

Now we can classify all conformal vectors in $U$.
First, we seek all conformal vectors with central charge $1/2$.
It is shown in \cite{M1} that there exists a one-to-one correspondence 
between the set of conformal vectors with central charge $c$ in $U$
and the set of idempotents with squared length $c/8$ in $U_2$.
So we should determine all idempotents with squared length $1/16$ in
$U_2$. 
Let $X=\alpha \w +\beta e +\gamma f +\delta f'$ be a conformal vector 
with central charge $1/2$.
Then by \cite{M1} we should solve the equations $(X/2)(1)(X/2)=(X/2)$ 
and $\la X,X\ra =1/16$.
By direct calculations, the solutions of $(\alpha,\beta,\gamma,\delta)$ 
are $(0,1,0,0)$, $(0,0,1,0)$ and $(0,0,0,1)$.
Therefore,

\begin{thm}\label{conformal}
  There are exactly three conformal vectors with central charge $1/2$
  in $U_2$, namely $e$, $f$ and $f'$.
\end{thm}

The rest of conformal vectors can be obtained in the following way.
As shown in \cite{M1}, conformal vectors and idempotents in the Griess 
algebra of $U$ are in one-to-one correspondence.
So we should seek all idempotents and their squared lengths in $U_2$. 
Since we have a set of basis $\{ \w, e, f,f'\}$ of $U_2$ and 
all multiplications and inner products are known, we can get 
them by direct calculations.
After some computations, we reach that the possible central charges are 
$1/2$, $81/70$, $58/35$, $4/5$ and $6/7$.
In the following, $(\alpha,\beta,\gamma,\delta )$ denotes 
$\alpha \w +\beta e +\gamma f +\delta f'$.
\vsb\\
(1)\ Central charge $1/2$: 
  $(0,1,0,0)$, 
  $(0,0,1,0)$, 
  $(0,0,0,1)$.
\vsb\\
(2)\ Central charge $81/70$: 
  $(1,-1,0,0)$, 
  $(1,0,-1,0)$, 
  $(1,0,0,-1)$.
\vsb\\
(3)\ Central charge $58/35$:
  $(1,0,0,0)$.
\vsb\\
(4)\ Central charge $4/5$: 
  $(14/9,-32/27,-32/27,-32/27)$, 
  $(-7/18,14/27,32/27,32/27)$, 
  $(-7/18,32/27,14/27,32/27)$, 
  $(-7/18,32/27,32/27,14/27)$.
\vsb\\
(5)\ Central charge $6/7$: 
  $(-5/9,32/27,32/27,32/27)$, 
  $(25/18,-14/27,-32/27,-32/27)$, 
  $(25/18,-32/27,-14/27,-32/27)$, 
  $(25/18,-32/27,-32/27,-14/27)$.

\subsection{Automorphisms} 

Let $V$ be any VOA and $e\in V$ a rational conformal vector with
central charge 1/2.
Then $e$ defines an involution $\tau_e$ of a VOA $V$, which is
so-called the Miyamoto involution (cf. \cite{M1}).
By Theorem \ref{conformal}, $U$ has three conformal vectors $e$, $f$
and $f'$.
Since $e^{\tau_f}\ne e$ nor $f$ and $f^{\tau_e}\ne f$ nor $e$, we must
have $e^{\tau_f}=f^{\tau_e}=f'$. 
Therefore, $\tau_{e}\tau_{f}\tau_{e}= \tau_{f^{\tau_e}}=
\tau_{e^{\tau_f}}= \tau_f\tau_e\tau_f$ and so $(\tau_e\tau_f)^3=1$.
It is clear that both $\tau_e$ and $\tau_f$ are non-trivial
involutions acting on $U$ and $\tau_e\ne \tau_f$.
Hence $\tau_e$ and $\tau_f$ generate $S_3$ in $\aut (U)$.
We prove that $\la \tau_e,\tau_f\ra = \aut (U)$.

\begin{thm}
  $\aut (U)=\la \tau_e,\tau_f\ra$.
\end{thm}

\pf
Let $g\in \aut (U)$.
Since $U$ is generated by $e$ and $f$, the action of $g$ on $U$ is
completely determined by its actions on $e$ and $f$.
By Theorem \ref{conformal}, the set of conformal vectors with central
charge $1/2$ in $U$ is $\{ e,f,f'\}$ so that we get an injection from
$\aut (U)$ to $S_3$. 
Since $\la \tau_e,\tau_f\ra$ acts on $\{ e,f,f'\}$ as $S_3$,
we obtain $\aut (U)=\la \tau_e,\tau_f\ra$.
\qed

\begin{rem}
  We note that both $\w^3$ and $\w^4$ are $S_3$-invariant so that the 
  orthogonal decomposition $\w =\w^3+\w^4$ is the {\it characteristic} 
  decomposition of $\w$ in $U$.
\end{rem}

Summarizing everything, we have already shown that $U$ is generated by 
two conformal vectors $e$ and $f$ whose inner product is 
$\la e,f\ra =13/2^{10}$ and its automorphism group is generated by 
two involutions $\tau_e$ and $\tau_f$ with $(\tau_e\tau_f)^3=1$.
Hence, we conclude that our VOA $U$ is the same as 
$\text{VA}(e,f)$ in \cite{M3} and gives a positive solution for
Theorem 5.6 (4) of \cite{M3}.

\begin{thm}\label{fixed-point}
  As a $\la \w^3\ra\tensor \la \w^4\ra$-module, 
  $U^{\la \tau_e,\tau_f\ra}=L(\fr{4}{5},0)\tensor L(\fr{6}{7},0)\oplus 
  L(\fr{4}{5},3)\tensor L(\fr{6}{7},5)$.
  It is a rational VOA.
\end{thm}

\pf
Since we may identify $U$ as $\text{VA}(e,f)$ in \cite{M3}, we can use 
the results obtained in \cite{M3}.
It is shown in \cite{M3} that $\la \w^3\ra \tensor \la \w^4\ra
=L(\fr{4}{5},0)\tensor L(\fr{6}{7},0)$ is a proper sub VOA of $U^{\la
\tau_e,\tau_f\ra}$. 
Since $U$ has both a $\Z_2$-grading and a $\Z_3$-grading, all
irreducible $L(\fr{4}{5},0)\tensor L(\fr{6}{7},0)$-submodules but
$L(\fr{4}{5},0)\tensor L(\fr{6}{7},0)$ and $L(\fr{4}{5},3)\tensor
L(\fr{6}{7},5)$ cannot be contained in $U^{\la \tau_e,\tau_f\ra}$.
Hence, $U^{\la \tau_e,\tau_f\ra}$ must be as stated.
The rationality of $U^{\la \tau_e,\tau_f\ra}$ will immediately follow
from results in \cite{LLY}.
\qed

\subsection{Fusion rules} 

Here we determine all fusion rules for irreducible $U$-modules.
We will denote the fusion product of irreducible $V$-modules $M^1$ and
$M^2$ by $M^1\fusion_V M^2$. 
Set $U=U^0\oplus U^1\oplus U^2$ with $U^0 = W(0)\tensor N(0)$, $U^1 =
W(\fr{2}{3})^+\tensor N(\fr{4}{3})^+$ and $U^2=W(\fr{2}{3})^-\tensor
N(\fr{4}{3})^-$. 
Recall the list of all irreducible $U$-modules shown in Theorem
\ref{irrs}.
We note that all of them are $\Z_3$-stable and each irreducible
$U$-module contains one and only one of the following irreducible
$U^0$-modules: 
$$
  W(h)\tensor N(k),\ h=0,\fr{2}{5},\ k=0,\fr{1}{7},\fr{5}{7}.
$$
Therefore, seen as $U^0$-modules, all irreducible $U$-modules have the
shape 
$$
\begin{array}{l}
  \ds
  U\fusion_{U^0} \l( W(h)\tensor N(k)\r)
  \vsb\\
  \ds
  = W(h)\tensor N(k)\oplus \{ U^1\fusion_{U^0} \l( W(h)\tensor N(k)\r) 
    \} \oplus \{ U^2\fusion_{U^0} \l( W(h)\tensor N(k)\r)\}
\end{array}
$$
with $h=0,\fr{2}{5}$ and $k=0,\fr{1}{7},\fr{5}{7}$.
Since $U\fusion_{U^0} \l( W(h)\tensor N(k)\r)$ denotes a $U^0$-module
in general, we denote an irreducible $U$-module of the shape
$U\fusion_{U^0} \l( W(h)\tensor N(k)\r)$ with $h=0,\fr{2}{5}$ and
$k=0,\fr{1}{7}, \fr{5}{7}$ by $\ind_{U^0}^U W(h)\tensor N(k)$ to
emphasize that it is a $U$-module.
Using this notation, the fusion products for irreducible $U$-modules
can be computed as follows: 

\begin{thm}
  All fusion rules for irreducible $U$-modules are given by the
  following formula:
  \begin{equation}\label{fusion_rule}
  \begin{array}{l}
    \ds
    \dim_\C
    \binom{\ind_{U^0}^U W(h_3)\tensor N(k_3)}{\ind_{U^0}^U W(h_1) 
    \tensor N(k_1)\q \ind_{U^0}^U W(h_2)\tensor N(k_2)}_U
    \vsv\\
    \ds \hspace{2cm}
    =\q 
    \dim_\C \binom{U\fusion_{U^0} \l( W(h_3)\tensor
    N(k_3)\r)}{W(h_1)\tensor N(k_1)\q W(h_2)\tensor N(k_2)}_{U^0},
  \end{array}
  \end{equation}
  where $h_1,h_2,h_3\in \{ 0,\fr{2}{5}\}$ and $k_1,k_2,k_3\in 
  \{ 0,\fr{1}{7}, \fr{5}{7}\}$.
\end{thm}

\pf
Since all irreducible $U$-modules are $\Z_3$-graded, the assertion
immediately follows from Lemma \ref{restriction} and Lemma
\ref{lifting}. 
\qed

\section{Application to the moonshine VOA} 

In this section, we work over the real number field $\R$.
We make it a rule to denote the complexification $\C \tensor_\R A$ of
a vector space $A$ over $\R$ by $\C A$.
Recall the construction of our VOA $U$ in Section 4.1.
In it, we only used a formula \eqref{character}, which was shown by
Goddard et al.{} by using a character formula in \cite{GKO}. 
Therefore, we can construct $U$ even if we work over $\R$.
To avoid confusions, we denote the real form of $U$ by $U_\R$.
We also note that the calculations on the Griess algebra of $U_\R$ in
Section 4.4 is still correct.

\subsection{Uniqueness theorem}

\begin{df}
  A VOA $V$ over $\R$ is referred to be {\it of moonshine type} if it 
  admits a weight space decomposition $V=\oplus_{n=0}^\infty V_n$ with
  $V_0=\R \vac$ and $V_1=0$ and it possesses a positive definite
  invariant bilinear form $\la\cd,\cd\ra$ such that $\la \vac,\vac\ra
  =1$. 
\end{df}

Assume that a VOA $V$ of moonshine type contains two distinct rational
conformal vectors $e$ and $f$ with central charge $1/2$.
In \cite{M3}, Miyamoto studied a vertex algebra $\text{VA}(e,f)$
generated by $e$ and $f$ in the case where their Miyamoto involutions
$\tau_e$ and $\tau_f$ generate $S_3$. 
In this subsection, we shall complete the classification of
$\text{VA}(e,f)$ in \cite{M3} in the case where the inner product $\la 
e,f\ra$ is $13/2^{10}$.

\begin{thm}\label{M3}
  \cite{M3}
  Under the assumption above, the inner product $\la e,f\ra$ is either 
  $1/2^6$ or $13/2^{10}$.
  When the inner product is equal to $13/2^{10}$, then a vertex
  algebra $\Vef$ generated by $e$ and $f$ forms a sub VOA in
  $V$.
  Denote by $\Vef^{(\tau_e\pm)}$ the eigen spaces for $\tau_e$ with
  eigenvalues $\pm 1$, respectively.
  The Griess algebra $\Vef_2$ is of dimension 4 and we can choose a
  basis 
  $\Vef_2^{(\tau_e+)} = \R \w^3 \perp \R \w^4 \perp \R v^0$ and
  $\Vef^{(\tau_e-)} = \R v^1$ such that $\w^3 +\w^4$ is the Virasoro
  vector of $\Vef$ and the multiplications and inner products in
  $\Vef_2$ are given as 
  $$
  \begin{array}{l}
    \ds
    \w^3_{(1)} \w^3=2\w^3,\q \w^3_{(1)} \w^4 =0,\q \w^3_{(1)} 
    v^0=\fr{2}{3}v^0,\q \w^3_{(1)} v^1=\fr{2}{3}v^1,\q 
    \w^4_{(1)} \w^4_{(1)}=2\w^4,
    \vsb\\
    \ds
    \w^4_{(1)} v^0=\fr{4}{3}v^0,\q
    \w^4_{(1)} v^1=\fr{4}{3}v^1,\q 
     v^0_{(1)} v^0_{(1)} = \fr{5}{6}\w^3 + \fr{14}{9}\w^4 -
    \fr{10}{9}v^0,\q v^0_{(1)}v^1=\fr{10}{9}v^1,
    \vsb\\
    \ds
    \la \w^3,\w^3\ra =\fr{2}{5},\q \la \w^4,\w^4\ra = \fr{3}{7},\q 
    \la v^0,v^0\ra =\fr{1}{2},\q \la v^1,v^1\ra =1.
  \end{array}
  $$
  The complexification $\C \Vef$ has a $\Z_3$-grading $\C \Vef =
  X^0\oplus X^1\oplus X^2$ and as $\C \mathrm{VA}(\w^3,\w^4)\simeq
  L(\fr{4}{5},0)\tensor L(\fr{6}{7},0)$-modules, they are isomorphic
  to one of the following:
  $$
  \begin{array}{ll}
  (i) & X^0 = \{ L(\fr{4}{5},0)\oplus L(\fr{4}{5},3)\} \tensor
    L(\fr{6}{7},0),\q 
    X^1 = L(\fr{4}{5},\fr{2}{3})^+ \tensor L(\fr{6}{7},\fr{4}{3}),
  \vsb\\
  & X^2 = L(\fr{4}{5},\fr{2}{3})^- \tensor L(\fr{6}{7},\fr{4}{3}); 
  \vsb\\
  (ii) & X^0 = L(\fr{4}{5},0)\tensor \{ L(\fr{6}{7},0)\oplus 
    L(\fr{6}{7},5)\},\q 
    X^1 = L(\fr{4}{5},\fr{2}{3})\tensor L(\fr{6}{7}, \fr{4}{3})^+, 
  \vsb\\
  & X^2=L(\fr{4}{5},\fr{2}{3})\tensor L(\fr{6}{7}, \fr{4}{3})^-;
  \vsb\\
  (iii) & X^0=L(\fr{4}{5},0)\tensor L(\fr{6}{7},0)\oplus
    L(\fr{4}{5},3) \tensor L(\fr{6}{7},5),\q
    X^1 = \{ L(\fr{4}{5},\fr{2}{3}) \tensor
    L(\fr{6}{7},\fr{4}{3})\}^+, 
  \vsb\\
  & X^2= \{ L(\fr{4}{5},\fr{2}{3}) \tensor
    L(\fr{6}{7},\fr{4}{3})\}^-;
  \vsb\\
  (iv) & X^0=\{ L(\fr{4}{5},0)\oplus L(\fr{4}{5},3)\} \tensor \{ 
    L(\fr{6}{7},0)\oplus L(\fr{6}{7},5)\},\q 
    X^1= L(\fr{4}{5},\fr{2}{3})^+ \tensor 
    L(\fr{6}{7},\fr{4}{3})^\pm, 
  \vsb\\
  & X^2= L(\fr{4}{5},\fr{2}{3})^-\tensor
    L(\fr{6}{7},\fr{4}{3})^\mp.
  \end{array}
  $$
  In the above, $M^-$ indicates a $\Z_2$-conjugate module of $M^+$. 
\end{thm}

We will prove the following.

\begin{thm}\label{completion}
  With reference to Theorem \ref{M3}, 
  only the case (iv) occurs.
  Therefore, $\C \Vef$ is isomorphic to $U=\C U_\R$ constructed in
  Section 4.
\end{thm}

\pf
First we prove (a).
The symmetric group $S_3=\la \tau_e,\tau_f\ra$ on three letters has
three irreducible representations $W_0=\C w^0$, $W_1=\C w^1$
and $W_2= \C w^2\oplus \C w^3$, where $W_0$ is a trivial module,
$\tau_e$ and $\tau_f$ act on $w^1$ as a scalar $-1$, and $\tau_e$ acts
on $w^2$ and $w^3$ as scalars respectively $1$ and $-1$. 
By quantum Galois theorem (cf.{} \cite{DM1} \cite{HMT}), 
we can decompose $\C \Vef$ as follows:
$$
  \C \Vef =\C \Vef^{\la \tau_e,\tau_f\ra} \tensor W_0 
  \bigoplus M_1\tensor W_1 \bigoplus M_2\tensor W_2,
$$
where $M_1$ and $M_2$ are inequivalent irreducible $\C\Vef^{\la
\tau_e,\tau_f\ra}$-modules. 
In the proof of Theorem \ref{M3} in \cite{M3}, Miyamoto found that
only the following two cases could be occur: $\C\Vef =\C
\mathrm{VA}(\w^3,\w^4)$ or $\C \Vef \supsetneq \C
\mathrm{VA}(\w^3,\w^4)$ and the former corresponds to the case
(i)-(iii) and the latter does the case (iv).
We assume that $\C\Vef =\C \mathrm{VA}(\w^3,\w^4)\simeq
L(\fr{4}{5},0)\tensor L(\fr{6}{7},0)$.
In this situation, $M^1$ as a $L(\fr{4}{5},0)\tensor
L(\fr{6}{7},0)$-module is isomorphic to $L(\fr{4}{5},3)\tensor
L(\fr{6}{7},0)$ in the case (i), $L(\fr{4}{5},0)\tensor
L(\fr{6}{7},5)$ in the case (ii) and $L(\fr{4}{5},3)\tensor
L(\fr{6}{7},5)$ in the case (iii) and $M^2$ as a
$L(\fr{4}{5},0)\tensor L(\fr{6}{7},0)$-module is isomorphic to
$L(\fr{4}{5},\fr{2}{3})\tensor L(\fr{6}{7},\fr{4}{3})$ in each case.
Therefore, $\C\Vef^{(\tau_e-)}$ has the following shapes:
$$
  \C\Vef^{(\tau_e-)}=
  \begin{cases}
    \, L(\fr{4}{5},3)\tensor L(\fr{6}{7},0)\tensor w^1 \bigoplus 
    L(\fr{4}{5},\fr{2}{3})\tensor L(\fr{6}{7},\fr{4}{3}) \tensor 
    w^3 & \mbox{in the case (i),}
    \vsb\\
    \, L(\fr{4}{5},0)\tensor L(\fr{6}{7},5) \tensor w^1 \bigoplus 
    L(\fr{4}{5},\fr{2}{3})\tensor L(\fr{6}{7},\fr{4}{3}) \tensor 
    w^3 & \mbox{in the case (ii),}
    \vsb\\
    \, L(\fr{4}{5},3)\tensor L(\fr{6}{7},5) \tensor w^1 \bigoplus
    L(\fr{4}{5},\fr{2}{3})\tensor L(\fr{6}{7},\fr{4}{3}) \tensor 
    w^3 & \mbox{in the case (iii).}
  \end{cases}
$$
We show that $\dim \C \Vef^{(\tau_e-)}_3=3$.
Since $\C \Vef^{(\tau_e-)}_2=\C v^1$ and $v^1$ is a highest weight
vector with highest weight $(\fr{2}{3},\fr{4}{3})$, $\w^3_{(0)}v^1$
and $\w^4_{(0)}v^1$ are linearly independent vectors in $\C
\Vef^{(\tau_e-)}_3$. We claim that $\{ \w^3_{(0)}v^1, \w^4_{(0)}v^1, 
v^0_{(0)}v^1\} $ is a set of linearly independent vectors in $\C
\Vef^{(\tau_e-)}_3$.
Set $x^1=\w^3_{(0)}v^1$, $x^2=\w^4_{(0)}v^1$ and $x^3=v^0_{(0)}v^1$.
Using the commutator formula $[a_{(m)},b_{(n)}]=\sum_{i\geq 0}
\binom{m}{i} (a_{(i)}b)_{(m+n-i)}$, invariant property $\la a_{(m)}b^1,
b^2\ra = \la b^1, a_{(-n+2)}b^2\ra$ for $a\in \C\Vef_2$ and an identity
$(a_{(0)}b)_{(m)}=[a_{(1)}, b_{(m-1)}]-(a_{(1)}b)_{(m-1)}$, we can
calculate all $\la x^i,x^j\ra$, $1\leq i,j\leq 3$ and it is a routine
work to check that $\det (\la x^i,x^j\ra)_{1\leq i,j\leq 3}\ne 0$.
Since $\Vef= \Vef^{(\tau_e+)}\perp \Vef^{(\tau_e-)}$, the
non-singularity of a matrix $(\la x^i,x^j\ra)_{1\leq i,j\leq 3}$
implies that $x^1$, $x^2$ and $x^3$ are linearly independent. 
Therefore, $\dim \C \Vef^{(\tau_e-)}_3=3$.
One can also see that 
$$
  v^0_{(0)}v^1 -\fr{5}{9}(\w^3_{(0)}+\w^4_{(0)})v^1
$$
is a non-zero highest weight vector for $L(\fr{4}{5},0)\tensor
L(\fr{6}{7},0)$ with highest weight $(3,0)$.
Thus, the possibility of $\C \Vef$ is only the case (i).
We next show that $\dim \C \Vef^{(\tau_e-)}_5=12$.
Set
$$
\begin{array}{l}
  y^1= \w^3_{(-2)} v^1,\q 
  y^2= \w^3_{(-1)} \w^3_{(0)} v^1,\q 
  y^3= \w^3_{(-1)} \w^4_{(0)} v^1,\q 
  y^4= \w^3_{(0)} \w^3_{(0)} \w^4_{(0)} v^1,
  \vsb\\
  y^5= \w^3_{(0)} \w^4_{(-1)} v^1,\q
  y^6= \w^3_{(0)} \w^4_{(0)} \w^4_{(0)} v^1,\q
  y^7= \w^4_{(-2)} v^1,\q
  y^8= \w^4_{(-1)} \w^4_{(0)} v^1,
  \vsb\\
  y^9= \w^3_{(-1)} \l( v^0_{(0)}-\fr{5}{9} \w^3_{(0)}-\fr{5}{9} 
    \w^4_{(0)} \r) v^1,\q 
  y^{10}= \w^3_{(0)}\w^3_{(0)} \l( v^0_{(0)}-\fr{5}{9}
  \w^3_{(0)}-\fr{5}{9} \w^4_{(0)} \r) v^1,\q 
  \vsb\\
  y^{11}= \w^4_{(-1)} \l( v^0_{(0)}-\fr{5}{9} \w^3_{(0)}-\fr{5}{9} 
    \w^4_{(0)} \r) v^1,\q 
  y^{12}= v^0_{(-2)}v^1.
\end{array}
$$
By a similar method used in computations of $\la x^i,x^j\ra$, we can
calculate all $\la y^i,y^j\ra$, $1\leq i,j\leq 12$, based on the
informations of the Griess algebra of $\Vef$ and it is also a routine
work to show that $\det (\la y^i,y^j\ra)_{1\leq i,j\leq 12}\ne 0$.
Therefore, $y^i$, $1\leq i\leq 12$, are linearly independent vectors
in $\C \Vef^{(\tau_e-)}_5$.
On the other hand, the dimension of the weight 5 subspace of the case
(i) is 11, which is a contradiction.
Therefore, $\C \Vef^{\la \tau_e,\tau_f\ra}\supsetneq
\C\mathrm{VA}(\w^3,\w^4)$, which leads the desired conclusion.
Since $\C\Vef$ and $\C U_\R$ have unique VOA-structures, $\C\Vef\simeq
\C U_\R =U$. 
\qed

\subsection{Embedding into the moonshine VOA}

Let $V_\R^\nat$ be the moonshine VOA \cite{FLM} over $\R$.
It is well-known that the full automorphism group of the moonshine VOA 
is the Monster $\M$, the largest sporadic finite simple group
(cf. \cite{FLM}).
Since $V_\R^\nat$ is (of course) a VOA of moonshine type, its weight
two subspace forms an commutative algebra, called the monstrous Griess
algebra. 
As shown \cite{C} and \cite{M1}, there is a one-to-one correspondence 
between the 2A-involutions of the Monster and conformal vectors with
central charge $1/2$ in $(V_\R^\nat)_2$.

There is a pair $(e,f)$ of conformal vectors with central charge
$1/2$ in $V_\R^\nat$ such that $\tau_e\tau_f$ defines a 3A-triality
in $\M$.
It is shown in \cite{C} that the inner product $\la e,f\ra$ of such a
pair is equal to $13/2^{10}$.
Therefore, the complexification of the moonshine VOA $\C V_\R^\nat$
contains a sub VOA isomorphic to $U$ by Theorem \ref{completion}.

\begin{thm}\label{3A}
  There exists a sub VOA isomorphic to $U$ in the complexificated
  moonshine VOA $\C V_\R^\nat$.
  Therefore, $\C V_\R^\nat$ contains both the 3-state Potts model
  $L(\fr{4}{5},0)\oplus L(\fr{4}{5},3)$ and the tricritical 3-state
  Potts model $L(\fr{6}{7},0)\oplus L(\fr{6}{7},5)$ and we can define 
  the 3A-triality in the Monster by the $\Z_3$-symmetries of the
  fusion algebras for these models.
\end{thm}


\begin{center}
  \gs{Acknowledgments}
\end{center}
\begin{quote}
  The authors would like to thank Professor Masahiko Miyamoto for
  stimulating discussions. 
  They also thank Professor Ching  Hung Lam for his kindness and
  informing them of his joint work with Professor Hiromichi Yamada. 
\end{quote}

\small
\setlength{\baselineskip}{12pt}

\end{document}